\documentclass{article}
\usepackage[leqno,intlimits]{amsmath}
\usepackage{amssymb}
\usepackage{amsthm}
\usepackage{hhline}
\usepackage{hyperref}
\input{pstricks}
\input{multido}

\newtheorem{tm}{Theorem}[section]
\newtheorem{lm}[tm]{Lemma}
\newtheorem{pr}[tm]{Proposition}
\newtheorem{cor}[tm]{Corollary}

\numberwithin{equation}{section}

\newcommand\kon{r} 
\newcommand\mom{m} 
\newcommand*{\un}[1]{\underline{#1}}
\newcommand*{\Zb}{\mathbb Z}
\newcommand*{\Rb}{\mathbb R}
\newcommand*{\om}{\omega}
\newcommand*{\ze}{\zeta}
\newcommand*{\de}{\delta}
\newcommand*{\la}{\lambda}
\newcommand*{\ba}{\begin{aligned}}
\newcommand*{\ea}{\end{aligned}}
\newcommand*{\be}{\begin{equation}}
\newcommand*{\ee}{\end{equation}}
\newcommand*{\e}[1]{\text{\rm e}^{#1}}
\newcommand*{\vr}{\varrho}
\newcommand*{\Ev}{{\bf E}}
\newcommand*{\Pv}{{\bf P}}
\newcommand*{\Vv}{{\text{\bf Var}}}
\newcommand*{\Oc}{\mathcal O}
\newcommand*{\lc}{\lceil}
\newcommand*{\rc}{\rceil}
\newcommand*{\lf}{\lfloor}
\newcommand*{\rf}{\rfloor}
\newcommand*{\di}{\,\text{\rm d}}
\newcommand*{\al}{\alpha}
\newcommand*{\ga}{\gamma}
\newcommand*{\wt}{\widetilde}
\newcommand*{\wih}{\widehat}
\newcommand*{\hop}{\bigskip\noindent}
\newcommand*{\hip}{\smallskip\noindent}
\newcommand*{\Ac}{\mathcal A}
\newcommand*{\Bc}{\mathcal B}

\providecommand{\abs}[1]{\left\vert#1\right\vert}

\begin{document}

\title{Order of current variance  and diffusivity  in the asymmetric simple exclusion process}
\author{M\'arton Bal\'azs\thanks{MTA-BME Stochastics Research Group}, Timo Sepp\"al\"ainen\thanks{University of Wisconsin-Madison\newline
M. Bal\'azs was partially supported by the Hungarian Scientific Research Fund (OTKA) grants T037685, K60708, TS49835, and F67729,  the Bolyai Scholarship of the Hungarian Academy of Sciences, and National Science Foundation grant DMS-0503650. Part of this work was completed while Bal\'azs was a visiting assistant professor at University of Wisconsin-Madison.\newline
T.\ Sepp\"al\"ainen was partially supported by National Science Foundation grants DMS-0402231 and DMS-0701091.}}
\maketitle
\begin{abstract}
We prove that the variance of the current across a characteristic 
is of order $t^{2/3}$ in a stationary  asymmetric simple exclusion process,
and that the diffusivity has order $t^{1/3}$.  The proof proceeds 
via couplings to 
show the corresponding moment bounds for a second class particle. 
\end{abstract}


\hop

\section{Introduction}

The  asymmetric
simple exclusion process (ASEP) is a Markov process that describes
the motion of particles  on the one-dimensional 
integer lattice $\Zb$, subject to the exclusion
interaction that allows at most one particle at each site.  
Particles in the process
jump one step to the right with rate $p$ and
one step to the left with rate $q=1-p$,
and we assume $0\le q< p\le 1$.  Particles attempt jumps
independently of each other, but any 
 attempt to jump onto an already occupied site is suppressed.
In Section \ref{sc:prore}
  below we give a rigorous construction of ASEP 
in terms of Poisson clocks that govern the jump attempts. 
This process is among the interacting particle systems introduced
in Spitzer's seminal paper \cite{spi}.    
We refer the reader to Liggett's
monographs  \cite{ips, stochi} for coverage of most of the 
 work on ASEP up to  the late 1990's.

In 1994 Ferrari and Fontes \cite{se}  proved
a central limit theorem 
for the net particle current    seen by an observer moving 
at a fixed speed $v$.  This quantity that we denote by
$J^{(v)}(t)$ is the number of particles that pass the
observer from left to right minus the number that 
pass from right to left during time interval  $(0,t]$. 
   The particle process is assumed 
to be stationary with Bernoulli-distributed occupation variables
at some average density $\vr\in(0,1)$. 
The result is the weak limit of the diffusively
rescaled and centered current: 
\[
\lim_{t\to\infty} \frac{J^{(v)}(t)-\Ev[J^{(v)}(t)]}{t^{1/2}}
\;=\; \chi_v.\]
 The limit $\chi_v$ is 
 a centered Gaussian random variable  with variance 
\[\sigma^2=\vr(1-\vr)\lvert (p-q)(1-2\vr)-v\rvert.\] 

The interesting phenomenon is 
 the vanishing of the variance at the characteristic speed
$v=V^\vr\equiv (p-q)(1-2\vr)$. As we explain below,
 $V^\vr$ is the
speed at which perturbations travel in the system, both at the  
microscopic particle level (in the  
expected sense as given in \eqref{eq:Qmean} below) 
and at the macroscopic p.d.e.\ level.  

Physical reasoning \cite{vBKS85} implied that the correct order of the 
fluctuations of the current across the characteristic
 should be $t^{1/3}$.  This would 
explain the degenerate limit under $t^{1/2}$ normalization. 
These ``$t^{1/3}$ fluctuations'' remained elusive throughout the 
1990's.  

The seminal papers of Baik, Deift and Johansson
\cite{bdj}  and Johansson \cite{1/3}  gave the first rigorous 
proofs of such fluctuations.   The correct order 
 was verified to be $t^{1/3}$, and the limiting fluctuations
were found to obey  Tracy-Widom distributions from
 random matrix theory. 
The first paper dealt with
the last-passage version of Hammersley's process, and 
the second with the last-passage version 
of the {\sl totally asymmetric} simple exclusion process (TASEP). 
 Total asymmetry means here that particles are
allowed to jump only in one direction at a constant rate,
so this is the case $p=1,q=0$. 

These papers did not study stationary particle processes, but
instead processes started from special jam-type deterministic 
initial conditions. For TASEP this means that initially
all sites to the left of the origin are occupied and all
sites to the right empty. 
With such initial conditions
the processes could be represented by 
 last-passage percolation models, a point that had been
exploited already in the past (among the seminal ones
were   \cite{aldo-diac},
\cite{sepp-ejp} and \cite{hkl}).  The actual
analysis was then performed entirely on 
combinatorial descriptions of the last-passage model.
  Later a last-passage representation was also found for a stationary
TASEP \cite{spohn}, and then the Tracy-Widom limit  proved 
for the current across the characteristic in that setting
\cite{ferspohn}. 

The early proofs 
of fluctuations relied on a counting argument that utilizes the 
Robinson-Schensted-Knuth correspondence for Young tableaux, and 
Gessel's formula that converts certain Schur function sums into
Toeplitz determinants. Later this step has been replaced by a more
direct connection between the last-passage model and a 
determinantal point process. 
 The fluctuation limits are then derived
by analyzing the asymptotics of the determinant in the appropriate 
scaling regime. 
Consequently, while a genuine breakthrough has been achieved,
the delicate steps of the proof restrict the reach of the results
in several ways.  In particular, the particles of the systems
are permitted to move in only one direction and admit only the
simplest type of jumps. 
 
In the present paper we give the first proof of the accurate
order of the fluctuations in systems that are only partially
asymmetric.  Namely, we 
show in the original setting of 
Ferrari-Fontes \cite{se} that the variance of the current
across the characteristic in 
 $(p,q)$ ASEP is of order $t^{2/3}$. 
Our arguments are entirely probabilistic and utilize couplings
of several processes and bounds on second class particles.
Informally speaking, {\sl second class particles} are 
perturbations in the system
that do not disturb the motion of the regular particles
but are influenced by the ambient system.  Precise 
definition must wait for  the construction
of the coupled processes in Section \ref{sc:prore}.    
Presently it appears  that there is no way to apply the 
combinatorial-analytic approach pioneered in \cite{bdj} 
to ASEP because there is no  last-passage model where the 
analysis could begin. 

We take a key insight from recent work of 
 Cator and Groeneboom \cite{cuberoot} and from our joint work
with Cator in  \cite{third}: this is the idea of coupling
processes whose densities  differ by $\mathcal O(t^{-1/3})$ in order
to bound the motion of a second class particle, whose fluctuations in turn
are linked to the variance of the current.   The couplings
we utilize go back to 
 Ferrari, Kipnis and Saada \cite{fks} who introduced them to
study the microscopic locations of shocks in the particle system. 

Fluctuation results for asymmetric exclusion processes  
have also been stated in terms of  a quantity called the 
{\sl diffusivity} $D(t)$.   One way to view the link between
current variance and diffusivity involves the
 second class particle:  the variance of
the current is the expected absolute deviation of the second
class particle, while $tD(t)$ is the variance of the second
class particle.  For ASEP we also obtain the correct 
order $t^{1/3}$ for the diffusivity.   In Section \ref{sc:prore}
we state the main result
Theorem \ref{tm:main} which  
is a moment bound for the second class particle. Bounds 
for current variance and diffusivity appear as Corollaries 
\ref{tm:curr} and \ref{tm:diff}. 

There is also work on the diffusivity with resolvent
methods.  The resolvent approach 
can  handle more general jump kernels than the 
nearest-neighbor type (the modifier ``simple'' in ASEP 
refers to the restriction to nearest-neighbor jumps). 
But so far this approach has not yielded optimal bounds. 
Results in one and two dimensions have  appeared in 
\cite{lqsy} and \cite{yau}.  

The most recent work utilizing the resolvent approach
is by  Quastel and Valk\'o \cite{quava, quava2}. 
In \cite{quava}   they show that, for any two finite-range 
exclusion processes with nonzero-mean jump distributions,
 the ratio of Laplace transforms of $tD(t)$
is uniformly bounded. 
In \cite{quava2}  this comparison theorem is paired up
with our Corollary \ref{tm:diff} to get the correct order
of the Laplace transform of $tD(t)$ for all these exclusion
processes, and with an additional argument also the correct
pointwise upper bound $D(t)\le Ct^{1/3}$. 

A few more words about the broader context.   We see this paper
as an opening for a treatment of several other models, such as 
zero-range and bricklayer processes.    For these systems
 diffusive current fluctuations off the characteristic
 were established by Bal\'azs in
\cite{fluct}.  

The $t^{1/3}$ fluctuations with Tracy-Widom
limits are universal for some class of asymmetric systems whose
precise characterization is not clear at the moment.  There is also
another universality class among asymmetric systems in one 
dimension, one where fluctuations occur on the scale $t^{1/4}$ 
and limits are Gaussian processes related to
fractional Brownian motion with Hurst parameter $H=1/4$.  Such
results appear in the papers \cite{dgl},  \cite{flucha}, and \cite{raprwre}.

\smallskip

{\sl Notation.}  We adhere to the usual notation for particle systems,
except that we underline notions that pertain to the entire lattice:
so $\un\eta(t)=\{\eta_i(t)\}_{i\in\Zb}$ denotes the occupation variables
of an exclusion process, and $\un\mu=\mu^{\otimes\Zb}$ is the 
i.i.d.\ product measure with marginal $\mu$. 
In general $\widetilde{X}$ denotes a centered random variable:
$\widetilde{X}=X-\Ev X$.  $C$ and $C_i$ for $i=1,2,3,\dotsc$
denote positive constants that can depend on $\vr$ and $p$ and can change
from line to line.  

\smallskip

{\sl Acknowledgement.} We thank an anonymous referee for a thorough
reading of the manuscript and insightful comments.  

\section{The exclusion process and the results}\label{sc:prore}


\noindent{\bf Construction of the process and second class particles}

\hip
The asymmetric simple exclusion process (ASEP) is a Markov 
process on the state space  $\Omega=\{0,\,1\}^\Zb$. Given a state 
$\un\om=\{\om_i\}_{i\in\Zb}\in\Omega$, the following jumps can happen 
independently at different sites:
\begin{align}
(\om_i,\,\om_{i+1})\longrightarrow(\om_i-1,\,\om_{i+1}+1)&\text{ with rate }
p\om_i(1-\om_{i+1}),\label{eq:up}\\
(\om_i,\,\om_{i+1})\longrightarrow(\om_i+1,\,\om_{i+1}-1)&\text{ with rate }
q(1-\om_i)\om_{i+1}.\label{eq:dn}
\end{align}
We assume $0\leq q=1-p<p\leq 1$. The special case $p=1$ is called
TASEP, or the totally asymmetric simple exclusion process.

We interpret the process as representing unlabeled particles
that execute independent nearest-neighbor random walks
on $\Zb$, subject to the exclusion interaction that suppresses 
 attempts to jump to an already occupied site.   
 $p$ is the rate of a particle 
to jump to the right and $q$ is the rate to jump to the left. 
The value $\om_i(t)=0$ means that  site $i$ is vacant at time $t$,
and  $\om_i(t)=1$  that  site $i$ is occupied at time $t$.
 The state of the entire process at time $t$ is 
then  $\un\om(t)$. 

A rigorous construction of this process is done by giving each site
$i$ two Poisson processes on the time line $[0,\infty)$: a rate
$p$ process $N_{i\to i+1}$ and a rate $q$ process  $N_{i\to i-1}$. 
The processes $\{ N_{i\to i+1},  N_{i\to i-1} :i\in\Zb\}$
are mutually independent, and also independent of the initial 
configuration $\un\om(0)$.   The rule of evolution is that 
when  $N_{i\to i+1}$ jumps, a particle is moved from $i$ to $i+1$ 
if $i$ is occupied and $i+1$ is vacant.  And similarly with 
 $N_{i\to i-1}$.  Thus the rates \eqref{eq:up}--\eqref{eq:dn}
are realized. 

Let $\mu_\vr$ denote the measure $\mu_\vr\{1\}=\vr=1-\mu_\vr\{0\}$ 
on the set $\{0,1\}$, and let $\un\mu_\vr=\mu_\vr^{\otimes\Zb}$
be the i.i.d.~Bernoulli product measure with
mean density $\vr$ on $\Omega$. 
It is known that the measures 
$\{\un\mu_\vr\,:\,0\leq\vr\leq1\}$ are the extreme points of 
the convex and weakly compact set of  
invariant distributions for the process
that are also invariant under spatial translations.

It is convenient to embed the exclusion process in a height
process that represents a wall of adjacent columns of bricks.
On top of each interval $[i,i+1]$ sits a  column of bricks with
 height $h_i\in\Zb$.  The entire height configuration
is $\un h=\{h_i\}_{i\in\Zb}$,  restricted to satisfy
\be 0\leq h_{i-1}-h_i\leq 1 \quad\text{ for each $i$} \label{eq:h}
\ee
so that the wall slopes downward to the right.
Let the Poisson processes govern the evolution of the heights:
 when $N_{i\to i+1}$ jumps add a brick
on top of the column on $[i,i+1]$, 
and when $N_{i+1\to i}$ jumps remove a brick
from the column on $[i,i+1]$.  But suppress every step
that  leads to a violation of \eqref{eq:h}. 
 
Given an initial particle configuration $\un\om(0)$, 
define an initial height configuration by
\[
h_i(0)=\begin{cases}
\sum_{j=i+1}^0\om_i(0)&\text{for }i<0,\\
0&\text{for }i=0,\\
-\sum_{j=1}^i\om_i(0)&\text{for }i>0.
\end{cases}\]
Let the heights evolve, and define
\[
\om_i(t)=h_{i-1}(t)-h_{i}(t).
\]
Then this process $\un\om(t)$ is exactly the ASEP constructed earlier,
and the height increment $h_i(t)-h_i(0)$ is the net particle 
current across the bond $(i,i+1)$. 

The Poisson construction reveals its power when it is used to 
run simultaneously several processes started from different 
initial states.  This is called the {\sl basic coupling}. The
first observation is that this coupling preserves monotonicity
among both particle and height configurations. 
Ordering is defined sitewise: for particle configurations
 $\un\eta\leq\un\om$ means that  $\eta_i\leq\om_i$  for each $i\in\Zb$,
and similarly for height configurations  
 $\un g\leq\un h$ if  $g_i\leq h_i$    for each $i\in\Zb$.
The basic  coupling has this property, called attractivity: 
\[
\un\eta(0)\leq\un\om(0)\Longrightarrow\un\eta(t)\leq\un\om(t)\quad\text{and}
\quad\un g(0)\leq\un h(0)\Longrightarrow\un g(t)\leq\un h(t)
\]
for all $t>0$.

We use the following terminology:  if we have 
two coupled exclusion processes such that  $\un\eta(t)\leq \un\om(t)$,
then the {\sl $\om-\eta$ second class particles} are the particles
that occupy sites $i$ at which $\om_i(t)-\eta_i(t)=1$.
The joint process $(\un\eta(\cdot),\un\om(\cdot))$ can be
constructed from a two-class process:
(i) The first class particles 
 $\un\eta$ obey the ASEP dynamics as described earlier. 
(ii) The second class particles $d_i=\om_i-\eta_i$ 
also obey the Poisson clocks when they can, but they are
not allowed to jump on sites occupied by first class particles,
and when a first class particle jumps on a site occupied by  a
second class particle, they swap sites. 

Let $\un\de_i\in\Omega$ denote a configuration that has only a single
particle  at 
site $i$. If $\un\eta\in\Omega$ is such that 
$\eta_0=0$, we can legitimately define
 $\un\eta^+=\un\eta+\un\de_0$. In this situation 
we say that there is  a single second class particle 
between $\un\eta^+$ and $\un\eta$ 
at site 0. Since the basic coupling conserves the single 
second class particle, there is always a site $Q(t)$ such that 
\begin{equation}
\un\eta^+(t)=\un\eta(t)+\un\de_{Q(t)}.\label{eq:delt}
\end{equation}
 $Q(t)$ is the position of the second class particle at time $t$, 
which performs a nearest neighbor walk, influenced by the ambient process 
$\un\eta(\cdot)$. 

It is convenient to also 
have the notion of a second class \emph{anti}particle at position 
$Q_a(t)$ in a process $\un\om(t)$.  This means  that 
$Q_a(t)$ is the location of the single discrepancy 
between two processes $\un\om(t)$ and  $\un\om^-(t)$ that are
 started so that 
$\un\om^-(0)=\un\om(0)-\un\de_i$ where $i=Q_a(0)$.
A moment's reflection reveals that in fact in the basic
coupling of ASEP a second class
particle and an antiparticle are the same thing.  But in the proofs
we will couple more than two processes and this extra flexibility will
be convenient. 


\hop
{\bf Current fluctuations and diffusivity}

\hip
Let $[x]$ denote the first integer from $x$ towards the origin, 
in other words  $[x]=\lf x\rf$ (floor) when $x\geq0$
 and $[x]=\lc x\rc$ (ceiling) 
when $x<0$. 
For a speed 
value $V\in\Rb$  define
\be
J^{(V)}(t)=h_{[Vt]}(t),\label{eq:jvdef}
\ee
the height of the column over interval $[\,[Vt],[Vt]+1]$ at time $t$.
Due to the normalization  $h_0(0)=0$, 
  $J^{(V)}(t)$ is the total net particle current seen by an observer moving
at speed $V$ during time interval $[0,t]$. Or more concretely, 
 $J^{(V)}(t)=J^{(V)}_+(t)-J^{(V)}_-(t)$ where 
$J^{(V)}_+(t)$ is the number of particles that
began in $(-\infty,0]$ at time $0$ but lie in $[\,[Vt]+1,\infty)$ 
at time $t$, and $J^{(V)}_-(t)$ is the number of particles that
began in $[1,\infty)$ at time $0$ but lie in $(-\infty, [Vt]\,]$ 
at time $t$. 
One can compute the density $\vr$ equilibrium expectation 
\[
\Ev J^{(V)}(t)= t(p-q)\vr(1-\vr)-\vr[Vt]
\]
by writing  a martingale for 
 $h_0(t)$ and then adding in
$h_{[Vt]}(t)-h_0(t)$ which counts particles 
between sites $0$ and $[Vt]$.  

Our results are based on an interplay between
currents and second class particles.  One key fact  is the
next  connection. 

\begin{pr}\label{pr:se}
Let $\un\om(\cdot)$ be an ASEP started 
from its stationary  Ber\-no\-ul\-li distribution  
$\un\mu_\vr$. Condition the origin to be empty and let 
$Q(\cdot)$ be a second class particle that starts at $Q(0)=0$.
Alternately condition the origin to be   occupied
and start a second class  antiparticle
$Q_a(\cdot)$ at the origin.  Either situation can be used 
to compute the variance of the current of the stationary process
for any $V\in\Rb$: 
\be
\begin{split}
\Vv(J^{(V)}(t))&=\vr(1-\vr)\Ev\bigl(\,|\,[Vt]-Q(t)|\,
\big\vert\,\om_0(0)=0\bigr)\\
 &=\vr(1-\vr)\Ev\bigl(\,|\,[Vt]-Q_a(t)|\,
\big\vert\,\om_0(0)=1\bigr) \end{split}
\label{eq:varjse}
\ee
\end{pr}
We shall also 
have occasion to use the following identity (true under the 
assumptions of the proposition above): 
\be
\Ev\bigl(Q(t)\,\bigr|\,\om_0(0)=0\bigr)=
\Ev\bigl(Q_a(t)\,\bigr|\,\om_0(0)=1\bigr)=t(p-q)(1-2\vr).
\label{eq:Qmean}\ee
Equations  \eqref{eq:varjse} and \eqref{eq:Qmean} have been derived and 
utilized earlier for TASEP,  \eqref{eq:varjse}  by Ferrari and Fontes \cite{se}
and  \eqref{eq:Qmean}  by Pr\"ahofer and Spohn \cite{spohn}. 
In article  \cite{varj2nd} identities \eqref{eq:varjse} and \eqref{eq:Qmean} 
are proved not only for ASEP but also for 
 zero-range and bricklayer processes.

The interesting current fluctuations 
occur at the \emph{characteristic speed} $V^\vr=(p-q)(1-2\vr)$.
From \eqref{eq:Qmean} we see that this is the average speed
of the second class particle.  To explain the 
significance of $V^\vr$ from another perspective,
we digress briefly to discuss a large scale,
 deterministic view of exclusion dynamics. 

  On large space and time
scales  the macroscopic particle density $\vr(t,x)$ of 
ASEP obeys the conservation law 
\be\vr_t+f(\vr)_x=0 \label{eq:conslaw}\ee 
 with flux $f(\vr)=(p-q)\vr(1-\vr)$.
This notion is made precise through a law of large numbers
called a {\sl hydrodynamic limit}. One considers a 
sequence of processes $\un\om^{(n)}(t)$ such 
that, for each $n$, the initial occupation variables 
are independent with means $\Ev\om^{(n)}_i(0)=\vr_0(i/n)$
for a given continuous function $\vr_0$. Then we
have the limit 
\be
\lim_{n\to\infty} \frac1n\sum_{na\le i\le nb}\om^{(n)}_i(nt)
=\int_a^b \vr(t,x)\,\di x
\quad\text{in probability} 
\label{eq:hydrolim}\ee
for all  intervals $[a,b]$.  The limiting density 
$\vr(t,x)$ is the entropy solution of \eqref{eq:conslaw}
with initial condition $\vr(0,x)=\vr_0(x)$.   
 Away from shocks
the {\sl characteristics}  of the partial differential equation 
\eqref{eq:conslaw} 
are solutions of the ordinary differential equation 
$\dot{x}=f'(\vr(t,x))$.  At constant density $\vr$ the 
characteristic speed is $f'(\vr)=V^\vr$.

  Limit
\eqref{eq:hydrolim} has been proved over time by a number of authors 
under various hypotheses.   Article
\cite{hl} was the first to achieve a satisfactory
level of generality for a broad class of asymmetric 
systems, also in higher  dimensions. 
For more information we refer  to  
\cite{cl} for 
 hydrodynamic limits
of interacting particle systems and 
to  \cite{evan} for 
 the theory
of scalar conservation laws. 

Let us return to the stationary ASEP with Bernoulli-$\vr$
occupation variables at each fixed time. 
A basic object for understanding space-time correlations
 is the  {\sl two point function}
$S(i,t)=\Ev[\om_i(t)\om_0(0)]-\vr^2$. 
It can  be written as the
transition probability of the second class particle:
\be\begin{split}
S(i,t) &= \vr\bigl( \Ev[\om_i(t)\,\vert\,\om_0(0)=1]-
\Ev[\om_i(t)] \,\bigr)\\
&=\vr(1-\vr)\bigl( \Ev[\om_i(t)\,\vert\,\om_0(0)=1]-
\Ev[\om_i(t)\,\vert\,\om_0(0)=0] \,\bigr)\\
&=\vr(1-\vr) \Ev[\,\om_i^+(t) -\om_i(t)
\,\vert\,\om^+_0(0)=1,\,\om_0(0)=0\,] \\
&= \vr(1-\vr) \Pv(Q(t)=i\,\vert\, \om_0(0)=0).
\end{split}\label{eq:SQcomp}\ee
On the third line above we coupled $\un\om(\cdot)$ with
a process $\un\om^+(\cdot)$ whose initial configurations
satisfy  $\un\om^+(0)=\un\om(0)+\un\delta_0$. In other
words,  the 
second class particle resides initially at the origin.
Note that use of the term ``transition probability'' 
is not meant to
suggest that $Q(t)$ is a Markov process.   

From \eqref{eq:SQcomp} and \eqref{eq:Qmean} follow the properties 
\[
 \sum_{i\in\Zb} S(i,t) = \vr(1-\vr)
\quad\text{and}\quad 
 \sum_{i\in\Zb}i S(i,t) = \vr(1-\vr) V^\vr t.
\]
 The  {\sl diffusivity}
is by definition a normalized second moment of the two-point
function: 
\be
D(t)=\frac1{t\vr(1-\vr)} \sum_{i\in\Zb} 
(i-V^\vr t)^2 S(i,t).
\label{eq:defDt}\ee
Consequently  the  diffusivity can also be expressed in terms
of the variance of the second class particle: 
\[
D(t)=t^{-1}\Vv(Q(t)\,|\,\om_0(0)=0).
\]
If the second class particle behaved like a random walk
(diffusively) then its variance would be of order $t$,
and thereby $D(t)$ of constant order.  The interesting
point is that the ASEP second class particle is 
{\sl superdiffusive} with variance of order $t^{4/3}$. 

We can now state the  main theorem,
 a moment bound for the second class 
particle in the setting where
 all the other occupation variables start off
Bernoulli-distributed. 

\begin{tm}\label{tm:main}
For ASEP $\un\om(\cdot)$ with rates  $0\leq q=1-p<p\leq 1$, started 
in $\un\mu_\vr$ distribution with any $\vr\in(0,1)$, there
exist constants $0<t_0, C<\infty$ such that for real $1\le \mom<3$
and  $t\ge t_0$ 
\be
C^{-1}\;\le\; \Ev\biggl\{\,
\Bigl\lvert\,\frac{Q(t)-V^\vr t}{t^{2/3}}\,\Bigr\rvert^\mom\;
\bigg\vert\,\om_0(0)=0\,\biggr\} \;\le\; C. 
\label{eq:mainbd}\ee
\end{tm}

 Let us
 simplify notation to $J^\vr(t)=J^{(V^\vr)}(t)$ for the
current across the characteristic. 
Via \eqref{eq:varjse}, 
taking $\mom=1$ gives the variance of this current.

\begin{cor}\label{tm:curr}
 With assumptions as in Theorem \ref{tm:main}, for $t\ge t_0$  
 the variance
of the current across the characteristic satisfies 
\[
C^{-1}t^{2/3} \le {\Vv(J^\vr(t))} \le C t^{2/3}.
\]
\end{cor}

Taking $\mom=2$ 
 identifies the order of the diffusivity. 

\begin{cor}\label{tm:diff}  With assumptions as in Theorem \ref{tm:main}, 
for $t\ge t_0$
\[
C^{-1}t^{1/3} \le D(t) \le C t^{1/3}.
\]
\end{cor}

The upper bound of \eqref{eq:mainbd} is not valid 
for all $t>0$. For small $t$ the variance of $Q(t)$ 
is of order $t$ because the second class particle
 is likely to have experienced at most
one  jump. 

The rest of the paper is devoted to the proof of
Theorem \ref{tm:main}, first the upper
bound and then the lower bound.  

\section{Upper bound}
\label{sc:ub}
Proof of the upper bound utilizes couplings of several 
processes and a tagged second class particle introduced
 by Ferrari, Kipnis and 
Saada \cite{fks}.  We refer the reader to the exposition of
these ideas in Section 
III.2 in Liggett's second monograph \cite{stochi}.
Before deriving the upper bound we introduce some 
preliminaries on the couplings. 


\hop
{\bf Couplings}

\hip
%
We start by describing an initial distribution
on three classes of particles and with labels attached to particles
of the two lower classes.  

Fix densities $0\leq\la<\vr<1$. Define the measure $\mu$ on $\{0,\,1\}\times
\{0,\,1\}$ by
\be
\mu(0,0)=1-\vr,\quad\mu(0,1)=\vr-\la,\quad\mu(1,0)=0,\quad
\text{and } \ \mu(1,1)=\la.
\label{eq:mupdef}
\ee
The first marginal of $\mu$ is $\mu_\la$ and the second is $\mu_\vr$. 
Let $\un\mu=\mu^{\otimes\Zb}$ be 
the i.i.d.\ product measure on $[\{0,\,1\}\times
\{0,\,1\}]^\Zb$ with marginals $\mu$. (Our notation 
for $\un\mu$ differs from that in \cite{stochi} where the 
second variable of $\un\mu$ is the 
process of second class particles.) 
We condition this measure  on having an
 $\om-\eta$ second class particle at the origin, and define
the measure  $\bar{\un\mu}$ by
\[
\bar{\un\mu}(\di\un\eta,\,\di\un\om)=
\un\mu(\di\un\eta,\,\di\un\om\,|\,\eta_0=0,\,\om_0=1).
\]
Then $\bar{\un\mu}$-a.s.\ we have the following conditions for $\un\eta$ and $\un\om$:
\begin{itemize}
\item $\un\eta\leq\un\om$,
\item there is a second class particle between $\un\eta$ and $\un\om$ at the 
origin, and there are infinitely many of them on both sides of the origin.
\end{itemize}
Label the  $\om-\eta$ second class particles with integers
in an increasing fashion from left to 
right, giving label 0 to the one initially at the origin.  
$X_k(0)$ denotes
the position of second class particle with label $k$, so that
\be\dotsm <X_{-1}(0)<X_0(0)=0<X_{1}(0)<\dotsm  
\label{eq:Xlabels}\ee

Given $(\un\eta, \un\om)$ chosen from distribution $\bar{\un\mu}$,
 define a third 
configuration $\un\ze$ as follows. 
 Let $\ze_i=\om_i$ whenever 
$\eta_i=\om_i$. At all other sites $i$ we have an
$\om-\eta$ second class particle: 
$0=\eta_i<\om_i=1$. If this second class particle has label $n$, then let
\be
\ze_i=\ze_{X_n(0)}=\left\{
\ba
\eta_i=0,&&\text{with probability }\frac{1}{1+(p/q)^n},\\
\om_i=1,&&\text{with probability }\frac{(p/q)^n}{1+(p/q)^n},
\ea
\right.\label{eq:zedef}
\ee
where these choices are independent from each other and everything else. 
For 
the case $p=1$ the probabilities are defined as the limits
 of the above formulas 
as $p=1-q\nearrow1$.

Let $\bar{\un\mu}^*$ denote the 
resulting joint distribution of $(\un\eta, \un\ze, \un\om)$. It is clear 
that   $\bar{\un\mu}^*$-a.s.\ we have  $\un\eta\leq\un\ze\leq\un\om$.

Let these configurations evolve
from the initial distribution $\bar{\un\mu}^*$
 with common Poisson clocks. 
This basic coupling preserves the ordering   $\un\eta(t)\leq\un\ze(t)
\leq\un\om(t)$ a.s.\ for all times $t$. 
The effect of the coupling is that the $\ze-\eta$ particles 
have priority over the  $\om-\ze$ particles. 

Consider the evolution of the  the $\om-\eta$ second class particles
in this coupling.  They start off labeled with integers as described
above \eqref{eq:Xlabels}. We let each particle keep its integer
label and let $X_k(t)$ denote the position of particle labeled $k$
at time $t$.     Since only
nearest-neighbor jumps are permitted, the ordering $X_k(t)<X_{k+1}(t)$ 
is preserved.  Additionally, particle $k$ carries a $\ze$-mark
$\zeta_{X_k(t)}(t)\in\{0,1\}$, initially assigned by \eqref{eq:zedef}.
$\zeta_{X_k(t)}(t)=0$ indicates that in the further subdivision
into second and third class particles, at time $t$ position
$X_k(t)$ is occupied by a third class particle. 
These $\ze$-marks are not retained but exchanged. 
 Whenever $X_k$ and $X_{k+1}$ are neighbors and carry different
$\ze$-marks $(1,0)$ or $(0,1)$, then $(1,0)$ becomes $(0,1)$ with
rate $p$ and  $(0,1)$ becomes $(1,0)$ with
rate $q$.    This is just a restatement of the effect of the 
basic coupling when a  $\ze-\eta$ particle is next to
an   $\om-\ze$ particle.  

The $\om-\eta$ second class particle that starts at the origin
is of special importance to us, so we set 
$X(t)=X_0(t)$.  

\begin{pr} Couple the processes as described above with initial
distribution  $\bar{\un\mu}^*$. 
Then for any fixed time $t\in[0,\infty)$ 
the marks $\{\zeta_{X_k(t)}(t): k\in\Zb\}$ are
 independent
of the entire process $\{(\un\eta(s), \un\om(s)):s\geq 0\}$. 
Furthermore, the distribution of marks is constant in time:
 given the evolution $(\un\eta(\cdot),\un\om(\cdot))$,
for all $0\leq t<\infty$ and disjoint sets $I_0$ and $I_1$ 
of integer labels, 
\begin{align*}
&\Pv \bigl\{ \text{$\zeta_{X_k(t)}(t)=0$ for $k\in I_0$
and $\zeta_{X_n(t)}(t)=1$ for $n\in I_1$ } \big\vert\, 
\un\eta(\cdot),\un\om(\cdot) \bigr\}  \\
&\qquad =\prod_{k\in I_0} \frac{1}{1+(p/q)^k} 
\; \cdot \; \prod_{n\in I_1} \frac{(p/q)^n}{1+(p/q)^n}.
\end{align*} 
\label{pr:fks}\end{pr}

This statement is hinted at 
in a remark after Proposition III.2.13 in \cite{stochi}. 
We summarize the argument for the reader's convenience.  
Observe first that we can construct the process of 
$(\un\eta,\un\om)$-particles
and $\ze$-labels with two independent collections of Poisson clocks, one that
governs the process $(\un\eta(\cdot),\un\om(\cdot))$, and 
another one that governs the exchanges of $\ze$-marks 
 among the $\om-\eta$ second class particles. 
This construction is not the same as the basic coupling described above,
but it leads to the same process because the infinitesimal
rates are the same. 

Consequently we can first construct the 
$(\un\eta(\cdot),\un\om(\cdot))$ evolution for all time, and then 
 superimpose on it the evolution of the $\ze$-marks.  The
claim is that, given the entire evolution
 $(\un\eta(\cdot),\un\om(\cdot))$,
the mark configuration $\{\zeta_{X_k(t)}(t):k\in\Zb\}$ on
the $\om-\eta$ second class particles  has
 distribution \eqref{eq:zedef}.  

The key point is reversibility.  Let us abbreviate 
temporarily $u_k(t)=\zeta_{X_k(t)}(t)$ for the mark of
second class particle $X_k$ at time $t$.  If the 
process $\un{u}(\cdot)$ 
obeyed the $(p,q)$-ASEP dynamics, the product distribution
with marginals \eqref{eq:zedef} would be one of the well-known 
blocking measures that are reversible for $\un{u}(\cdot)$.  
(This would in fact be the case if $\la=0$ and $\vr=1$, for
then the second class particles would occupy all sites
and never move: $X_k(t)=k$ for all $0\leq t<\infty$.)  
The 
result now follows from this observation: the evolution
$(\un\eta(\cdot),\un\om(\cdot))$ can be thought of as 
a dynamical ``environment'' for the mark process
that  admits or prohibits certain exchanges at different times:
$X_k$ and $X_{k+1}$ cannot exchange marks unless they 
occupy adjacent sites, and the time intervals 
during which this happens are determined by 
$(\un\eta(\cdot),\un\om(\cdot))$.  However, imposing such an
environment on the process does not change the reversibility
of the measure.  

This last point can be checked rigorously
by first letting only finitely many marks
  $U_K(t)=\{u_k(t): -K\leq k\leq K\}$ evolve while the remaining marks
are frozen.  Then we are talking about a finite state
 Markov chain.  Given the evolution  
$(\un\eta(\cdot),\un\om(\cdot))$, we can partition the time
axis  $0=t_0<t_1<t_2<\dotsm <t_i
\nearrow\infty$  
so that the generator of the chain $U_K(\cdot)$ does
not change during $(t_{i-1},t_i)$.  Since detailed balance
is not violated by symmetrically prohibiting certain jumps, reversibility
and hence invariance is preserved during each interval $(t_{i-1},t_i)$,
and thereby for all time.
Letting $K\to\infty$ extends the invariance to the 
infinite mark process.   This argument establishes Proposition
\ref{pr:fks}.

\hop

The purpose of the construction is to confine certain 
special particles to be introduced shortly. 
Notice that a.s.\ there are only finitely many 
positive labels $n>0$ for which the 
first event of \eqref{eq:zedef} happens.
 Hence we have a rightmost position $R(0)=\sup\{ i:  
\eta_i=\ze_i<\om_i\}$. This is the initial position of the rightmost 
second class particle between $\un\ze$ and $\un\om$. Similarly, the second 
event in \eqref{eq:zedef} happens a.s.\ for 
only finitely many negative $n<0$. 
Let  $L(0)=\inf\{i: \eta_i<\ze_i=\om_i\}$ be the leftmost position 
where this happens, the initial position 
of the leftmost second class particle between $\un\eta$ and $\un\ze$. 
Define the events
\[
\Ac=\{0\leq R(0)\}\qquad\text{and}\qquad\Bc=\{L(0)\leq0\leq R(0)\}.
\]
These events only depend on the marks
\eqref{eq:zedef}, and hence are independent of 
 $(\un\eta, \un\om)$. $\Ac$ always has positive probability, 
while $\Bc$ has positive probability if $p<1$.

These positions evolve in the coupling, so at time $t$  we let
 $R(t)$ be  the position of the rightmost second class 
particle between $\un\ze(t)$ and $\un\om(t)$, 
and similarly  $L(t)$ 
is the position of the leftmost second class particle between 
$\un\eta(t)$ and $\un\ze(t)$.  $L(t)$ and $R(t)$ are always among
the $\om-\eta$ second class particles $\{X_k(t)\}$.  
Let $n_L(t)$ and $n_R(t)$ be the random labels such that 
$L(t)=X_{n_L(t)}(t)$ and $R(t)=X_{n_R(t)}(t)$.  These labels
are functions of the marks:  
\be
n_L(t)=\inf\{k :\zeta_{X_k(t)}(t)=1\}
\quad\text{and}\quad
n_R(t)=\sup\{k :\zeta_{X_k(t)}(t)=0\}.
\label{eq:defnR} \ee
By the above proposition $n_L(t)$ and $n_R(t)$ are 
independent of  $(\un\eta(t),\,\un\om(t))$ and their
distribution does  not change with time. 

To complete the construction we define two more 
initial configurations by
\[
\un\eta^+=\un\eta+\un\de_0\qquad\text{and}\qquad\un\om^-=\un\om-\un\de_0 
\]
where $\un\de_0$ denotes a configuration with a single particle 
at the origin, as explained  in Section \ref{sc:prore}. 
The basic coupling then applies to all five 
processes $\un\eta(\cdot)$, $\un\eta^+(\cdot)$, $\un\ze(\cdot)$, $\un\om^-(\cdot)$, 
$\un\om(\cdot)$ with the above initial data. The orderings
\be
\un\eta(t)\leq\un\ze(t)\leq\un\om(t),\qquad\un\eta(t)\leq\un\om^-(t)\leq\un\om(t),
\qquad\un\eta(t)\leq\un\eta^+(t)\leq\un\om(t)\label{eq:io}
\ee
are preserved by the evolution. There is always one difference between 
$\un\eta(\cdot)$ and $\un\eta^+(\cdot)$ which is thought of as a single second 
class particle on $\un\eta(\cdot)$ in the sense \eqref{eq:delt}. Denote its 
position by $Q(t)$ with $Q(0)=0$. Also, $\un\om(\cdot)$ and $\un\om^-(\cdot)$ 
always have a single second class particle between them, thought of as a second 
class antiparticle on $\un\om(\cdot)$ as described
 in Section \ref{sc:prore}. Its 
position is $Q_a(t)$, and $Q_a(0)=0$.  See
Figure \ref{hmmwhatanicepicture} for an illustration of the
initial configurations in this coupling.

\psset{unit=1pt}
\begin{figure}[ht]
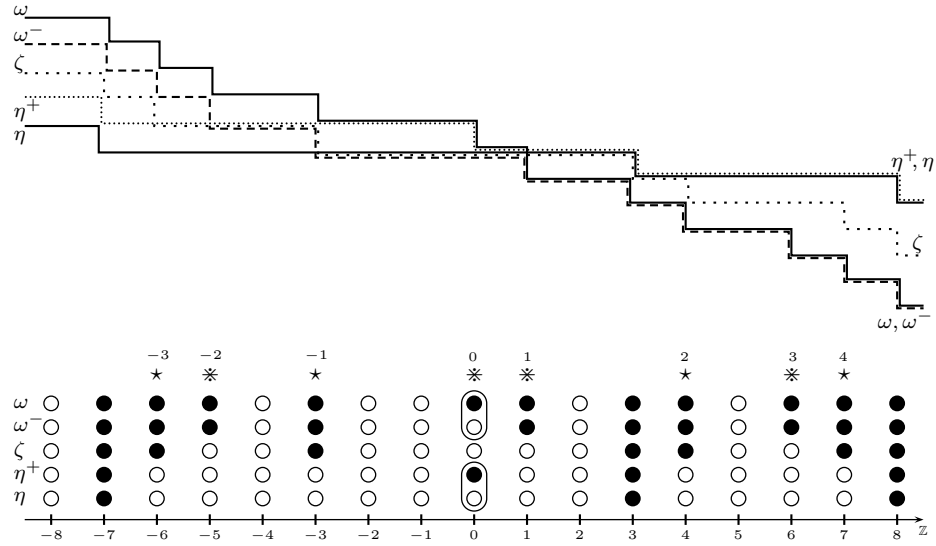

\begin{center}
\pspicture*(-15,0)(335,210)


\psline(-10,200)(22,200)(22,191)(41,191)(41,181)(61,181)(61,171)(101,171)(101,161)(161,161)(161,151)(180,151)(180,139)(219,139)(219,130)(240,130)(240,120)(280,120)(280,110)(301,110)(301,101)(321,101)(321,91)(330,91)

\psline[linestyle=dashed,dash=3pt 2pt](-10,190)(21,190)(21,180)(40,180)(40,170)(60,170)(60,158)(100,158)(100,147)(179,147)(179,138)(218,138)(218,129)(239,129)(239,119)(279,119)(279,109)(300,109)(300,100)(320,100)(320,90)(330,90)

\psline[linestyle=dotted,dotsep=1pt](-10,170)(19,170)(19,160)(160,160)(160,150)(222,150)(222,141)(321,141)(321,131)(330,131)

\psline(-10,159)(18,159)(18,149)(221,149)(221,140)(320,140)(320,130)(330,130)

\psline[linestyle=dashed,dash=1pt 3pt](-10,179)(20,179)(20,170)(39,170)(39,159)(101,159)(101,148)(220,148)(220,139)(241,139)(241,130)(300,130)(300,120)(320,120)(320,110)(330,110)

\rput[lb](-14,201){\footnotesize$\om$}
\rput[lb](-14,191){\footnotesize$\om^-$}
\rput[lb](-14,161){\footnotesize$\eta^+$}
\rput[lt](-14,157){\footnotesize$\eta$}
\rput[lb](-14,179){\footnotesize$\ze$}

\rput[rb](334,81){\footnotesize$\om,\om^-$}
\rput[rb](334,142){\footnotesize$\eta^+\!\!,\eta$}
\rput[rb](330,111){\footnotesize$\ze$}

\psline[linewidth=0.5pt]{->}(-10,10)(330,10)
\multido{\i=0+20,\n=-8+1}{17}{
\rput[t](\i,6){\tiny$\n$}
\psline(\i,8)(\i,12)
}
\rput(330,6){\tiny$\Zb$}

\rput[l](-14,18){\footnotesize$\eta$}
\pscircle[linewidth=0.5pt](0,18){3}
\qdisk(20,18){3}
\pscircle[linewidth=0.5pt](40,18){3}
\pscircle[linewidth=0.5pt](60,18){3}
\pscircle[linewidth=0.5pt](80,18){3}
\pscircle[linewidth=0.5pt](100,18){3}
\pscircle[linewidth=0.5pt](120,18){3}
\pscircle[linewidth=0.5pt](140,18){3}
\pscircle[linewidth=0.5pt](160,18){3}
\pscircle[linewidth=0.5pt](180,18){3}
\pscircle[linewidth=0.5pt](200,18){3}
\qdisk(220,18){3}
\pscircle[linewidth=0.5pt](240,18){3}
\pscircle[linewidth=0.5pt](260,18){3}
\pscircle[linewidth=0.5pt](280,18){3}
\pscircle[linewidth=0.5pt](300,18){3}
\qdisk(320,18){3}

\rput[l](-14,28){\footnotesize$\eta^+$}
\pscircle[linewidth=0.5pt](0,27){3}
\qdisk(20,27){3}
\pscircle[linewidth=0.5pt](40,27){3}
\pscircle[linewidth=0.5pt](60,27){3}
\pscircle[linewidth=0.5pt](80,27){3}
\pscircle[linewidth=0.5pt](100,27){3}
\pscircle[linewidth=0.5pt](120,27){3}
\pscircle[linewidth=0.5pt](140,27){3}
\qdisk(160,27){3}
\pscircle[linewidth=0.5pt](180,27){3}
\pscircle[linewidth=0.5pt](200,27){3}
\qdisk(220,27){3}
\pscircle[linewidth=0.5pt](240,27){3}
\pscircle[linewidth=0.5pt](260,27){3}
\pscircle[linewidth=0.5pt](280,27){3}
\pscircle[linewidth=0.5pt](300,27){3}
\qdisk(320,27){3}

\rput[l](-14,36){\footnotesize$\ze$}
\pscircle[linewidth=0.5pt](0,36){3}
\qdisk(20,36){3}
\qdisk(40,36){3}
\pscircle[linewidth=0.5pt](60,36){3}
\pscircle[linewidth=0.5pt](80,36){3}
\qdisk(100,36){3}
\pscircle[linewidth=0.5pt](120,36){3}
\pscircle[linewidth=0.5pt](140,36){3}
\pscircle[linewidth=0.5pt](160,36){3}
\pscircle[linewidth=0.5pt](180,36){3}
\pscircle[linewidth=0.5pt](200,36){3}
\qdisk(220,36){3}
\qdisk(240,36){3}
\pscircle[linewidth=0.5pt](260,36){3}
\pscircle[linewidth=0.5pt](280,36){3}
\qdisk(300,36){3}
\qdisk(320,36){3}

\rput[l](-14,47){\footnotesize$\om^-$}
\pscircle[linewidth=0.5pt](0,45){3}
\qdisk(20,45){3}
\qdisk(40,45){3}
\qdisk(60,45){3}
\pscircle[linewidth=0.5pt](80,45){3}
\qdisk(100,45){3}
\pscircle[linewidth=0.5pt](120,45){3}
\pscircle[linewidth=0.5pt](140,45){3}
\pscircle[linewidth=0.5pt](160,45){3}
\qdisk(180,45){3}
\pscircle[linewidth=0.5pt](200,45){3}
\qdisk(220,45){3}
\qdisk(240,45){3}
\pscircle[linewidth=0.5pt](260,45){3}
\qdisk(280,45){3}
\qdisk(300,45){3}
\qdisk(320,45){3}

\rput[l](-14,54){\footnotesize$\om$}
\pscircle[linewidth=0.5pt](0,54){3}
\qdisk(20,54){3}
\qdisk(40,54){3}
\qdisk(60,54){3}
\pscircle[linewidth=0.5pt](80,54){3}
\qdisk(100,54){3}
\pscircle[linewidth=0.5pt](120,54){3}
\pscircle[linewidth=0.5pt](140,54){3}
\qdisk(160,54){3}
\qdisk(180,54){3}
\pscircle[linewidth=0.5pt](200,54){3}
\qdisk(220,54){3}
\qdisk(240,54){3}
\pscircle[linewidth=0.5pt](260,54){3}
\qdisk(280,54){3}
\qdisk(300,54){3}
\qdisk(320,54){3}

\psframe[linewidth=0.5pt,framearc=5](165,32)(155,13)
\psframe[linewidth=0.5pt,framearc=5](165,59)(155,40)

\rput(40,65){\small$\star$}
\rput(60,65){\footnotesize$\divideontimes$}
\rput(100,65){\small$\star$}
\rput(160,65){\footnotesize$\divideontimes$}
\rput(180,65){\footnotesize$\divideontimes$}
\rput(240,65){\small$\star$}
\rput(280,65){\footnotesize$\divideontimes$}
\rput(300,65){\small$\star$}

\rput[lb](36,70){\tiny$-3$}
\rput[lb](56,70){\tiny$-2$}
\rput[lb](96,70){\tiny$-1$}
\rput[lb](158,70){\tiny$0$}
\rput[lb](178.6,70){\tiny$1$}
\rput[lb](238,70){\tiny$2$}
\rput[lb](279,70){\tiny$3$}
\rput[lb](298,70){\tiny$4$}

\endpspicture
\end{center}
\caption{\small 
A realization of the initial state of the five coupled processes.
The top portion of the figure shows the height functions and
the bottom portion the corresponding
 particle configurations. The lattice $\Zb$
runs at the bottom.   Solid disks denote particles
 and  open circles denote vacant sites.  Above the particles, 
{\footnotesize$\divideontimes$} marks locations of  
 $\om-\ze$  particles 
(first case in \eqref{eq:zedef}) and {\small$\star$} 
marks locations of 
$\ze-\eta$ particles (second case in \eqref{eq:zedef}). 
Together these make up the $\om-\eta$ second 
class particles and the numbers above these symbols represent
the labeling \eqref{eq:Xlabels}. 
 Both events $\Ac$ and $\Bc$ happen in 
this picture. The particle-hole pairs inside the ovals
represent the initial locations
of $Q_a$ and $Q$.  
\label{hmmwhatanicepicture}}
\end{figure}

The bounds will be proved by controlling the particles $Q(t)$ 
and $Q_a(t)$.  The next lemma contains one of the key points.

\begin{lm} In the coupling of five processes we have these 
implications:
\be
\Ac\subset\{Q_a(t)\leq R(t)\text{ for all }t\geq0\}.\label{eq:defac}
\ee
and 
\be
\Bc\subset\{Q_a(t)\leq R(t)\text{ and }L(t)\leq Q(t)\text{ for all }t\geq0\}.
\label{eq:defbc}
\ee
\label{lm:AB}
\end{lm}

\begin{proof}  We prove \eqref{eq:defac} and leave the similar argument
for \eqref{eq:defbc} to the reader.

Under the event $\Ac$, $Q_a(0)=0\leq R(0)$ holds initially. That is, we have 
a $\om-\zeta$ second class particle at or to the right of
 the second class antiparticle  $Q_a$. 
We show 
that the coupling preserves the inequality $Q_a(t)\leq R(t)$.
We have two cases to consider that could potentially
allow  $Q_a> R$ to happen. 
\begin{itemize}
\item 
Suppose the ordering $\un\ze(s)\leq\un\om^-(s)\leq\un\om(s)$ holds
at some time $s$. 
Since the ordering is then preserved 
 for all later times $t\geq s$, $\ze_{Q_a(t)}(t)\leq\om^-_{Q_a(t)}(t)
<\om_{Q_a(t)}(t)$, which implies an $\om-\ze$ second class particle 
at location $Q_a(t)$. As $R(t)$ is the position of the rightmost such 
second class particle, $Q_a(t)\leq R(t)$ holds.
\item Consequently the only possibility for producing $Q_a>R$
is to have a jump in a situation of this type: 
 $Q_a(t)=i$, $R(t)=i+1$ for some site $i$, and $\om^-_i(t)=0$, 
$\ze_i(t)=\om_i(t)=1$, $\ze_{i+1}(t)=0$, $\om^-_{i+1}(t)=\om_{i+1}(t)=1$. 
In this case column $i$ of $\un\ze$ increases by one with rate $p$, or 
column $i$ of $\un\om^-$ decreases by one with rate $q$. Neither
one of these 
steps interchanges $Q_a$ and $R$. 
\end{itemize}
This proves \eqref{eq:defac}. 
\end{proof}

For applications of this lemma it is crucial that 
the events $\Ac$ and $\Bc$ depend only on the initial marks \eqref{eq:zedef}
which were independent of the initial configuration $(\un\eta,\,\un\om)$.
Consequently  $\Ac$ and $\Bc$ are independent of 
the joint evolution  $(\un\eta(\cdot),\,\un\om(\cdot),\,X(\cdot),\,Q(\cdot),
\,Q_a(\cdot))$. 
When we condition the processes on $\Ac$ or $\Bc$
we call  this construction a \emph{conditional coupling}.


\hop
{\bf Proof of the upper bound}

\hip
Throughout this section the probability measure $\Pv$ 
represents the five-process coupling constructed in the previous
section.  
We begin with the proof of the upper bound in Theorem \ref{tm:main}. 
$C$ and its 
variants $C_1, C_2, \dotsc$ 
 denote positive constants that possibly depend on 
$p$ and $\vr$ 
and whose values can change from line to line. 
We first prove a lower bound on 
the density of the $\om-\eta$ second class particles. 
For integers $j\in\Zb$ and 
$u>0$  let
\be
N_j(t)=\sum_{i=j+1}^{j+2u}(\om_i(t)-\eta_i(t)).\label{eq:ntdef}
\ee
\begin{lm}\label{lm:2ndldp}  Let $\la, \vr\in(0,1)$ and $d\geq 0$ be an
integer.  Then there are strictly positive finite
 constants $\ga=\ga(\vr)$, $C_1=C_1(\vr,d)$ and 
$C_2=C_2(\vr)$  such that the following holds:
if  $0<\vr-\la<\ga$, then 
 for all integers $j\in\Zb$ 
and  $u>0$  and any time $t\geq0$,
\[
\Pv\Bigl\{N_j(t)<u(\vr-\la)+d\Bigr\}\leq C_1\exp\{-C_2u(\vr-\la)^2\}.
\]
\end{lm}
\begin{proof}
For the moment, denote by $\un y(\cdot)$ a 
$(p,q)$ exclusion process such that 
$y_i(0)=\om_i(0)$ at all sites $i$ except for $i=0$, and $\un z(\cdot)$ is a 
$(p,q)$ exclusion process such that $z_i(0)=\eta_i(0)$ at all sites $i$ except 
for $i=0$. For $i=0$, we pick the pair $(z_0(0),\,y_0(0))$ in distribution 
$\mu$ \eqref{eq:mupdef}, independently of the configuration on other sites. 
Apply the basic coupling to ensure $\un\eta(t)\leq\un z(t)\leq\un y(t)
\leq\un\om(t)$ for all $t\geq0$ (notice that this holds initially). $\un y(t)$ 
and $\un z(t)$ are marginally time-stationary processes,
 hence we can omit the 
notation for their time dependence in our arguments. However, the pair 
$(\un z(t),\,\un y(t))$ is \emph{not} in product distribution for $t>0$. Define
\[
Y=\sum_{i=j+1}^{j+2u}y_i\qquad\text{and}\qquad Z=\sum_{i=j+1}^{j+2u}z_i,
\]
so that $N_j(t)\geq Y-Z$. Then for any $\al>0$,
\begin{multline*}
\ba
\Pv\Bigl\{N_j(t)<u(\vr-\la)+d\Bigr\}&\leq\Pv\Bigl\{\e{-\al(Y-Z)}>
\e{-\al u(\vr-\la)-d\al}\Bigr\}\\
&\leq\e{\al u(\vr-\la)+d\al}\cdot\Ev\bigl(\e{-\al Y}\e{\al Z}\bigr)\\
&\leq\e{\al u(\vr-\la)+d\al}\cdot\bigl[\Ev\bigl(\e{-2\al Y}\bigr)\bigr]^{1/2}
\cdot\bigl[\Ev\bigl(\e{2\al Z}\bigr)\bigr]^{1/2}\\
&=\e{\al u(\vr-\la)+d\al}\cdot\bigl[\Ev\bigl(\e{-2\al y_0}\bigr)\bigr]^u
\cdot\bigl[\Ev\bigl(\e{2\al z_0}\bigr)\bigr]^u\\
\ea\\
\ba
&=\exp\bigl\{\al u(\vr-\la)+d\al+u\log\bigl[1+\vr(\e{-2\al}-1)\bigr]
+u\log\bigl[1+\la(\e{2\al}-1)\bigr]\bigr\}\\
&\leq\exp\bigl\{-\al u(\vr-\la)+d\al+2\al^2u(\vr+\la-\vr^2-\la^2)+uC\al^3\bigr\}
\ea
\end{multline*}
Here we used the marginal $\un\mu_\la$ and $\un\mu_\vr$ product distributions 
of $\un z$ and $\un y$, while the last inequality comes from Taylor expansion 
w.r.t.\ $\al$. The $\Oc(\al^3)$ term is uniform over $\la<\vr$ for a fixed $\vr$. 
Next we pick
\[
\al=\frac{\vr-\la}{4(\vr+\la-\vr^2-\la^2)},
\]
this optimizes the $\vr$ and $\la$-dependent terms, and finishes the proof.
\end{proof}

We turn to develop the main estimate for the upper bound.  The
objective is to bound the deviation 
$\Pv\{Q_a(t)\geq u+[V^\vr t]\}$ with a suitable expression that
involves the moment $\Ev\lvert Q_a(t)-[V^\vr t]\rvert$.  This is 
reached in \eqref{eq:final} below and then the upper bound comes 
from an elementary integration step.  

Along the way we compare currents in two processes 
that we abbreviate as follows:
\[
\ba
J^\vr(t)&=J^{((p-q)(1-2\vr))}(t)\quad\text{for current in the 
$\un\om(\cdot)$  process, and}\\
J^{V^\vr,\la}(t)&=J^{((p-q)(1-2\vr))}(t)\quad\text{for current in the 
$\un\eta(\cdot)$  process.}
\ea
\]
Notice that both use the same speed $V^\vr=f'(\vr)=(p-q)(1-2\vr)$
for the observer.  As already defined in the Introduction, 
this is the 
characteristic speed of the $\un\om(\cdot)$ process with density $\vr$. 
Let $A=(\Pv\{\Ac\})^{-1}<\infty$ from the 
conditional coupling, and recall that a tilde centers a random variable.
\begin{lm}
 Suppose $\vr-\la<\ga$ with $\ga$
from Lemma \ref{lm:2ndldp}.  Then for positive integers $u$
and times $t\in[0,\infty)$, 
\be
\ba
&\Pv\{Q_a(t)\geq4u+[V^\vr t]\}\\
&\qquad \leq A\Pv\{\wt J^\vr(t)-\wt J^{V^\vr,\la}(t)\geq 
u(\vr-\la)-t(p-q)(\vr-\la)^2\}\\
&\qquad\qquad +C_1\exp\{-C_2u(\vr-\la)^2\}.
\ea\label{eq:qjprob}
\ee
\end{lm}
\begin{proof}  By the independence observed after the proof of 
Lemma \ref{lm:AB}, conditioning on $\Ac$ does not affect the 
probability of the event of interest. 
By \eqref{eq:defac},
\begin{align*}
&\Pv\{Q_a(t)\geq4u+[V^\vr t]\}=\Pv\{Q_a(t)\geq4u+[V^\vr t]\,|\,\Ac\}\\
&\quad \leq\Pv\{R(t)\geq4u+[V^\vr t]\,|\,\Ac\}
\leq A\Pv\{R(t)\geq4u+[V^\vr t]\}\\
&\quad \leq A\Pv\{X(t)\geq2u+[V^\vr t]\}
+A\Pv\{R(t)\geq4u+[V^\vr t],\ X(t)<2u+[V^\vr t]\}.
\end{align*}
Recall that $N_{[V^\vr t]}(t)$ of \eqref{eq:ntdef} counts the number of 
$\om-\eta$ second class particles at time $t$ in the interval 
$\{[V^\vr t]+1,\dots,[V^\vr t]+2u\}$.  
Since the second class particles stay 
ordered and $X(t)$ started at the origin, 
the event $\{X(t)\geq2u+[V^\vr t]\}$ 
implies that all these second class particles crossed the path
 $s\mapsto [V^\vr s]+1/2$ by time $t$. 
 Each such second class 
particle crossing increases $J^\vr(t)-J^{V^\vr,\la}(t)$ by one. Therefore
\begin{multline*}
\Pv\{X(t)\geq2u+[V^\vr t]\}\leq\Pv\{J^\vr(t)-J^{V^\vr,\la}(t)
\geq N_{[V^\vr t]}(t)\}\\
\leq\Pv\Bigl\{J^\vr(t)-J^{V^\vr,\la}(t)\geq u(\vr-\la)+3\Bigr\}+\Pv\Bigl
\{N_{[V^\vr t]}(t)<u(\vr-\la)+3\Bigr\}.
\end{multline*}
Combine the previous displays to get
\begin{align}
&\Pv\{Q_a(t)\geq4u+[V^\vr t]\}\leq A\Pv\Bigl\{J^\vr(t)-J^{V^\vr,\la}(t)
\geq u(\vr-\la)+3\Bigr\}\label{eq:temp7}\\
&\qquad+A\Pv\Bigl\{N_{[V^\vr t]}(t)<u(\vr-\la)+3\Bigr\}\label{eq:nupr}\\
&\qquad+A\Pv\{R(t)\geq4u+[V^\vr t],\ X(t)<2u+[V^\vr t]\}.\label{eq:rxpr}
\end{align}
To line \eqref{eq:nupr} apply Lemma \ref{lm:2ndldp} with $d=3$. 
The event in \eqref{eq:rxpr} can happen in two ways:
\begin{itemize}
\item The label $n_R(t)$  of $R(t)$ is larger than $u(\vr-\la)$. 
The probability of 
this is, according to \eqref{eq:zedef} and its time-invariance, bounded by
\[
\sum_{n=\lc u(\vr-\la)\rc}^\infty\frac{1}{1+(p/q)^n}\leq\sum_{n=\lc u(\vr-\la)\rc}^\infty(q/p)^n\leq\frac{p}{p-q}\cdot(q/p)^{u(\vr-\la)}.
\]
\item There are fewer than $u(\vr-\la)$ $\om-\ze$ second class 
particles in the discrete interval $\{[V^\vr t]+2u,\dots,[V^\vr t]+4u-1\}$. The 
probability of this event is bounded in Lemma \ref{lm:2ndldp}.
\end{itemize}
Hence line \eqref{eq:rxpr} is bounded by 
\[
C_1\exp\{-uC_2(\vr-\la)^2\}+{Ap}{(p-q)^{-1}}\exp\{-uC_3(\vr-\la)\}.
\]
By modifying the constants,  we see that 
lines  \eqref{eq:nupr} and \eqref{eq:rxpr} are bounded by the last 
exponential term in  \eqref{eq:qjprob}. 

If $\un\om$ and $\un\eta$ start from their respective 
$\un\mu_\vr$ and $\un\mu_\la$ equilibria, then we would have
\be
\begin{split}
\Ev(J^\vr(t))-\Ev(J^{V^\vr,\la}(t))
&=t\bigl( f(\vr)-f(\la)-f'(\vr)(\vr-\la)\bigr)\\ 
&=-\tfrac12 t f''(\vr)(\vr-\la)^2\\
&=t(p-q)(\vr-\la)^2, 
\end{split}\label{eq:H''point}\ee
where we utilized the precise form $f(\vr)=(p-q)\vr(1-\vr)$
of the flux and 
ignored the error coming from the integer part of $V^\vr t$. 
This last  error is at most one. Our processes are also perturbed initially 
at the origin by the conditioning in $\un\mu$, which gives an error 
not larger 
than 2 in the above quantity.  
The term $+3$ inside the probability on line \eqref{eq:temp7}
 makes up for these errors. 
\end{proof}

{\sl Remark.} We spelled out calculation \eqref{eq:H''point} 
because the nonvanishing of $f''$ is a crucial factor that 
produces the order $t^{2/3}$ for the fluctuations of the second
class particle.  In  the end the variable $u$ is of order $t^{2/3}$
and $\vr-\la$ of order $t^{-1/3}$.

\begin{lm}
There is a constant $C>0$ such that for any time $t\geq0$,
\[
\Ev(|X(t)-L(t)|)\leq C(\vr-\la)^{-1}\quad\text{and}\quad\Ev(|R(t)-X(t)|)
\leq C(\vr-\la)^{-1}.
\]
\end{lm}
\begin{proof}
By Proposition III.2.10 of Liggett \cite{stochi}, the joint distribution 
$\bar{\un\mu}\overline{S(t)}$ of the pair $(\un\eta(t),\,\un\om(t))$, 
as seen from $X(t)$, can also be obtained by conditioning the distribution 
$\un\mu S(t)$ of the pair $(\un\eta(t),\,\un\om(t))$ at time $t$ on having 
a second class particle at the origin. (The result is denoted by 
$\overline{\un\mu S(t)}$ in \cite{stochi}.) The measure $\un\mu S(t)$ is 
translation-invariant, and gives probability $\vr-\la$ for having a second 
class particle at the origin. If $X_k(t)$ denotes the position of the 
$k^\text{th}$ second class particle at time $t$, then by Theorem B47 in 
 \cite{stochi}, we have
\[
\ba
\Ev(X_{k+1}(t)-X_k(t))&=\Ev^{\bar{\un\mu}\overline{S(t)}}(X_{k+1}-X_k)=
\Ev^{\overline{\un\mu S(t)}}(X_{k+1}-X_k)\\
&=\Ev^{\un\mu S(t)}(X_{k+1}-X_k\,|\,\om_0-\eta_0=1)\\
&=(\vr-\la)^{-1}\Ev^{\un\mu S(t)}\bigl(
(X_{k+1}-X_k)\cdot{\bf 1}\{\om_0-\eta_0=1\}\bigr)\\
&=\frac{1}{\vr-\la}.
\ea
\]

Recall that $R(t)=X_{n_R(t)}(t)$ with $n_R(t)$ defined by
\eqref{eq:defnR}  and  $X(t)=X_0(t)$. 
 Use  the independence and the time-invariance of the 
 distribution of marks   from  
Proposition \ref{pr:fks}:
\[
\ba
\Ev(|R(t)-X(t)|)&=\Ev\Bigl(\;\sum_{k=-\infty}^{-1}(X_{k+1}(t)-X_k(t))\cdot{\bf1}\{n_R(t)
\leq k\}\Bigr)\\
&\qquad+\Ev\Bigl(\;\sum_{k=0}^{\infty}(X_{k+1}(t)-X_k(t))\cdot{\bf1}\{k<n_R(t)\}\Bigr)\\
&=\sum_{k=-\infty}^{-1}\Ev(X_{k+1}(t)-X_k(t))\cdot\Pv\{n_R(t)\leq k\}\\
&\qquad+\sum_{k=0}^{\infty}\Ev(X_{k+1}(t)-X_k(t))\cdot{\Pv}\{k<n_R(t)\}\\
&=\frac{1}{\vr-\la}\cdot\Ev(|n_R(t)|)=\frac{C}{\vr-\la}.
\ea
\]
Similar considerations 
give the result for $L(t)$.
%
%
\end{proof}
\begin{lm}\label{lm:vardif}
\[
\Vv(J^{V^\vr,\la}(t))\leq\la\Ev\bigl(|[V^\vr t]-Q_a(t)|\bigr)+2t(p-q)\la(\vr-\la)+C(\vr-\la)^{-1}
\]
\end{lm}
\begin{proof}
The variance $\Vv$
 in the statement is taken in the five-process coupling 
where the Bernoulli  distribution  $\un\mu_\la$
 is initially perturbed at the origin. 
 Denote by $\Vv^\la$  
variance in the stationary process with initial invariant distribution 
 $\un\mu_\la$. 
By the conditional variance formula
\be
\begin{split}
\Vv^\la(J^{V^\vr,\la}(t))&=\Ev^\la\Vv^\la(J^{V^\vr,\la}(t)\,|\,\eta_0(0))\\
&\qquad +\Vv^\la\Ev^\la(J^{V^\vr,\la}(t)\,|\,\eta_0(0))\\
&\geq(1-\la)\cdot\Vv^\la(J^{V^\vr,\la}(t)\,|\,\eta_0(0)=0)\\
&=(1-\la)\cdot\Vv(J^{V^\vr,\la}(t)).
\end{split}\label{eq:temp12}\ee
Apply Proposition \ref{pr:se} to the first term 
$\Vv^\la(J^{V^\vr,\la}(t))$.  The 
conditional expectations on the right-hand side of \eqref{eq:varjse} 
match exactly the marginals of the five-process coupling
constructed earlier in this Section, and so we come back to the present setting:
\begin{align*}
\Vv(J^{V^\vr,\la}(t))  
&\leq\frac1{1-\la}\Vv^\la(J^{V^\vr,\la}(t))
=\frac1{1-\la}\la(1-\la)\Ev\bigl(|[V^\vr t]-Q(t)|\bigr)\\
&\leq\la\Ev\bigl(|[V^\vr t]-Q_a(t)|\bigr)+\la\Ev\bigl(|Q(t)-Q_a(t)|\bigr).
\end{align*}
Write  the last term above as 
\be
\Ev\bigl(|Q(t)-Q_a(t)|\bigr)=\Ev(Q(t))-\Ev(Q_a(t))+2\Ev([Q(t)-Q_a(t)]^-).
\label{eq:qnegp}
\ee
 Proposition \ref{pr:se} gives 
\[
\Ev(Q(t))-\Ev(Q_a(t))=2t(p-q)(\vr-\la).
\]

 The last term in \eqref{eq:qnegp}
is treated separately for TASEP ($p=1$) and ASEP ($0<q=1-p<p<1$). 

Consider first the totally asymmetric case $p=1$. 
Initially $Q(0)=Q_a(0)=0$, and we show $Q(t)\geq Q_a(t)$ 
for all $t\geq 0$. The following situations cover all cases 
where  $Q$ and $Q_a$ could interchange positions:
\begin{itemize}
\item If $Q(t)=Q_a(t)=i$, then at site $i$ \eqref{eq:io} has the unique 
solution $\eta_i(t)=\om_i^-(t)=0$ and $\eta_i^+(t)=\om_i(t)=1$.
 A right step of $Q_a$ without $Q$ would require a brick on column $i$ of $\un\om$, but not of $\un\eta^+$. This is impossible by the basic coupling and $\eta^+_{i+1}(t)=\eta_{i+1}(t)\leq\om_{i+1}(t)$. A left step of $Q$ without $Q_a$ would require a brick on column $i-1$ of $\un\eta$, but not of $\un\om^-$. This is again impossible by $\eta_{i-1}(t)\leq\om_{i-1}(t)=\om^-_{i-1}(t)$.
\item If $Q_a(t)=i$ and $Q(t)=i+1$, then \eqref{eq:io} has the unique solution $\eta_i(t)=\eta_i^+(t)=\om^-_i(t)=0$,  $\om_i(t)=1$, 
$\eta_{i+1}(t)=0$, and  $\eta_{i+1}^+(t)=\om_{i+1}^-(t)=\om_{i+1}(t)=1$. 
None of the processes can have column $i$ grow
 in this situation, hence $Q_a$ and $Q$ cannot interchange positions.
\end{itemize}
We conclude that the second term on the right of \eqref{eq:qnegp} is zero 
in the totally asymmetric case.

For ASEP we use conditional coupling with event $\Bc$ 
of \eqref{eq:defbc}. Let $B=(\Pv\{\Bc\})^{-1}<\infty$.
\[
\ba
&\Ev([Q(t)-Q_a(t)]^-)=\Ev([Q(t)-Q_a(t)]^-\,|\,\Bc)
\leq\Ev([L(t)-R(t)]^-\,|\,\Bc)\\
&\quad\leq B\Ev([L(t)-R(t)]^-)\leq B\Ev(|R(t)-X(t)|)+B\Ev(|X(t)-L(t)|)\\
&\quad \leq C(\vr-\la)^{-1}
\ea
\]
by the previous lemma.

Collecting terms completes the proof of the lemma.
\end{proof}

{\sl Remark.} The constant $B=(\Pv\{\Bc\})^{-1}$ diverges as
$p\nearrow 1$ and hence the same is true of $C$ in the last term
of the bound in Lemma \ref{lm:vardif}.  This suggests that there
should be a more effective way to bound the variance of 
$J^{V^\vr,\la}(t)$ of the $\la$-system in terms of quantities
computed in the $\vr$-system. 

\begin{lm}
\[
\Vv(J^\vr(t)-J^{V^\vr,\la}(t))\leq2\Ev\bigl(|[V^\vr t]-Q_a(t)|\bigr)+4t(p-q)\la(\vr-\la)+C(\vr-\la)^{-1}.
\]
\end{lm}
\begin{proof}
Similarly to the previous proof, we write
\[
\Vv^\vr(J^\vr(t))\geq\vr\cdot\Vv(J^\vr(t)).
\]
Then we proceed by
\[
\ba
&\Vv(J^\vr(t)-J^{V^\vr,\la}(t))\leq2\Vv(J^\vr(t))+2\Vv(J^{V^\vr,\la}(t))\\
&\qquad \leq{2}{\vr^{-1}}\Vv^\vr(J^\vr(t))+2\Vv(J^{V^\vr,\la}(t))\\
&\qquad \leq2(1-\vr+\la)\Ev\bigl(|[V^\vr t]-Q_a(t)|\bigr)
+4t(p-q)\la(\vr-\la)+C(\vr-\la)^{-1}
\ea
\]
utilizing  the previous lemma.
\end{proof}

We come to the lemma that summarizes all the previous estimation. 

\begin{lm}
For any real $u\geq 1$ and time $t> 0$,
\begin{multline}
\Pv\{Q_a(t)\geq4u+[V^\vr t]\}\leq C_3\frac{t^2}{u^4}\cdot\Ev\bigl(|Q_a(t)-[V^\vr t]|\bigr)\\
+C_4\frac{t^2}{u^3}+C_5\frac{t^3}{u^5}
+C_1\exp\Bigl\{-C_2\frac{u^3}{t^2}\Bigr\}+\e{-2u}.\label{eq:final}
\end{multline}
\end{lm}
\begin{proof}  
Set $b=2(\gamma\wedge\vr)(p-q)$ where $\gamma$ is the constant from
Lemma \ref{lm:2ndldp}.  We proceed by cases. 

{\sl Case 1:}  $1\le u<bt$.  Suppose first $u$ is an integer
as it was in the proof of \eqref{eq:qjprob}.  Throughout
density $\vr$ has been fixed, and now we also fix
\[
\la=\vr-\frac{u}{2t(p-q)}.
\]
The constraint on $u$ guarantees that $\la>0$ and 
$\vr-\la<\ga$ which was required for \eqref{eq:qjprob}. The point
of this choice of $\la$ is to maximize
 the  lower bound inside the 
probability on the right-hand side of \eqref{eq:qjprob}.
So, continuing from \eqref{eq:qjprob} with 
 Chebyshev's inequality and the previous lemma, 
\begin{multline*}
\Pv\{Q_a(t)\geq4u+[V^\vr t]\}\\
\ba
&\leq A\Pv\Bigl\{\wt J^\vr(t)-\wt J^{V^\vr,\la}(t)\geq\frac{u^2}{4t(p-q)}\Bigr\}+
C_1\exp\Bigl\{-C_2\frac{u^3}{t^2}\Bigr\}\\
&\leq C_3\frac{t^2}{u^4}\Vv(J^\vr(t)-J^{V^\vr,\la}(t))+
C_1\exp\Bigl\{-C_2\frac{u^3}{t^2}\Bigr\}\\
&\leq C_3\frac{t^2}{u^4}\cdot\Ev\bigl(|Q_a(t)-[V^\vr t]|\bigr)+
C_4\frac{t^2}{u^3}+C_5\frac{t^3}{u^5}+C_1\exp\Bigl\{-C_2\frac{u^3}{t^2}\Bigr\}.
\ea
\end{multline*}
Extension  from integral $u$ to real $u$ is achieved by adjusting  
constants on the last line above. 

{\sl Case 2:}  $bt\le u< 2t$. Then $bu/2< bt$ and 
\[
\Pv\{Q_a(t)\geq4u+[V^\vr t]\}\leq\Pv\{Q_a(t)\geq4\cdot bu/2+[V^\vr t]\}.
\]
{\sl Case 1}  can be applied
with $u$ replaced by $bu/2$, at the price of adjusting
some constants with powers of $b$. 

{\sl Case 3:}  $u\ge 2t$.  Since  $Q_a(t)$ is 
bounded above by a rate one Poisson process, 
\[
\Pv\{Q_a(t)\geq4u+[V^\vr t]\}\leq\Pv\{Q_a(t)\geq3u\}\leq\e{-2u}
\]
for all $t$.  Combining the bounds from the cases proves the lemma.
\end{proof}

We are ready to complete the proof of 
the upper bound of Theorem \ref{tm:main}. 
Abbreviate \[ \Psi(t)=\Ev\bigl(\,|Q_a(t)-[V^\vr t]|\,\bigr).\]

First we need to complement \eqref{eq:final} with a matching
lower tail bound for $Q_a(t)$. This can be derived by 
arguments analogous to the ones we have pursued
throughout Section \ref{sc:ub}. 

Alternately, we can derive the lower tail bound by a
 particle-hole interchange followed by reflection  of the lattice. 
Define $\wih\om_i(t)=1-\om_i(t)$, an ASEP
 with density $\wih\vr=1-\vr$ and rightward jump rate $\wih p=1-p$. 
Its second class 
particle has position $\wih Q(t)=Q_a(t)$.  
Let $\wih{\un\om}^R(\cdot)$
be the process obtained from  $\wih{\un\om}(\cdot)$ through
a reflection about 
the origin.  Then $\wih{\un\om}^R(\cdot)$ is 
ASEP  with the original parameters $(p,q)$. The second class 
particle of this process 
is at position $\wih Q(t)^R=-\wih Q(t)=-Q_a(t)$, and the 
characteristic speed is 
$\wih V^{\vr, R}=(p-q)(1-2(1-\vr))=-V^\vr$. Since $[-V^\vr t]=-[V^\vr t]$,
the expectation 
\[
\Psi(t)=\Ev\bigl(|Q_a(t)-[V^\vr t]|\bigr)
=\Ev\bigl(|\wih Q(t)^R-[\wih V^{\vr, R}t]|\bigr)
\]
is the same for $\un\om(\cdot)$ and $\wih{\un\om}^R(\cdot)$. 
Consequently 
\[
\Pv\{Q_a(t)-[V^\vr t]\le -4u\}=\Pv\{\wih Q(t)^R-[V^{\vr,R} t]\ge 4u\}
\]
which is again bounded by the right-hand side of \eqref{eq:final}.

Introduce  a large constant  $4<\kon<\infty$. 
Consider $t\ge 1$ and  $u\ge \kon t^{2/3}$.  Then we can 
combine the upper and lower tail bounds whose common
right-hand side is given in \eqref{eq:final}, replace $4u$
by $u$ (we made sure $u/4\ge 1$),  and simplify a little to arrive at 
\be
\begin{split}
&\Pv\{\,\lvert Q_a(t)-[V^\vr t] \rvert \ge u\} \\
&\qquad \leq C_1\frac{t^2}{u^4} \Psi(t)
+C_2(\kon)\Bigl(\,\frac{t^2}{u^3} + 
\exp\Bigl\{-C_3(\kon)\frac{u}{t^{2/3}}\Bigr\} \,\Bigr).
\end{split}\label{eq:final2}\ee
Constants were renamed and their dependence on  $\kon$ expressed
in the notation  because now it is of importance that $C_1$ in front of 
$\Psi(t)$ does not depend on $\kon$. 
In particular, the exponentials in \eqref{eq:final} were 
combined via 
$e^{-2u}$ $\le$ $\exp(-2ut^{-2/3})$ and 
$\exp(-C_2u^3t^{-2})\le \exp(-C_2\kon^{2}ut^{-2/3})$.    

Let $1\le \mom<3$.  Integrate bound \eqref{eq:final2} over
$u\in [\kon t^{2/3},\infty)$:
\be\begin{split}
&\Ev\bigl(\,|Q_a(t)-[V^\vr t]|^\mom\,\bigr)\\
&\qquad \le \kon^\mom t^{2\mom/3} +\mom\int_{\kon t^{2/3}}^\infty 
\Pv\{\,\lvert Q_a(t)-[V^\vr t] \rvert \ge u\} u^{\mom-1} \,\di u\\
&\qquad \le C_1 \kon^{\mom-4} \Psi(t) t^{2+(2/3)(\mom-4)} \;+\;
C_4(\kon)t^{2\mom/3}. 
\end{split}\label{eq:final3}\ee 
Constants were renamed again and their dependence on $\mom$ ignored. 
To get the final bounds, take first $\mom=1$ 
in \eqref{eq:final3} to get 
\[
\Psi(t) \le C_1 \kon^{-3} \Psi(t)  \;+\; C_4(\kon)t^{2/3}. 
\]
Since $C_1$ is independent of $\kon$,
fixing  $\kon$ large enough gives $\Psi(t)\le C_5(\kon)t^{2/3}$.
Putting this bound back on the last line of \eqref{eq:final3}
then gives   
\[
\Ev\bigl(\,|Q_a(t)-[V^\vr t]|^\mom\,\bigr) \le C_6(\kon)t^{2\mom/3}
\]
for $1<\mom<3$. 
The upper bound of Theorem \ref{tm:main} has been proved.

\section{Lower bound}

The lower bound is proved by perturbing an initial equilibrium
on a segment of the lattice.  Again we begin with a description and 
some properties of the coupling. 


\hop
{\bf Perturbing a segment initially}

\hip
Recall again 
the characteristic speeds $V^\vr=(p-q)(1-2\vr)$ and $V^\la=(p-q)(1-2\la)$. 
We assume $\vr>\la$, hence $V^\vr<V^\la$. 
Throughout this section  $u>0$ denotes a fixed positive integer, and 
\[
n=[V^\la t]-[V^\vr t]+u.
\]
To begin  define an initial 
product distribution on two configurations $(\un\eta(0),\un\ze(0))$
 by describing the marginals
on each lattice site $i$: 
\[\begin{cases}
(\eta_i(0),\,\ze_i(0))\sim\mu\text{ of \eqref{eq:mupdef}}&\text{if }i<-n,\\
\eta_i(0)=0,\ \ze_i(0)\sim\mu_\vr&\text{if }i=-n,\\
\eta_i(0)=\ze_i(0)\sim\mu_\la&\text{if }-n<i\leq0,\\
(\eta_i(0),\,\ze_i(0))\sim\mu\text{ of \eqref{eq:mupdef}}&\text{if }i>0.
\end{cases} 
\]
Let  $Q^{(-n)}(t)$ be the position at time $t$ of 
 a second class particle in the process  $\un\eta(\cdot)$, 
initially   at site $Q^{(-n)}(0)=-n$.  Note that 
site $-n$ was set vacant for $\un\eta(0)$. 
As before $\un\eta(t)$ together with the particle at $Q^{(-n)}(t)$
make up  the process 
$\un\eta^+(t)$. Except for this perturbation $\un\eta(0)$ starts 
 in the 
Bernoulli $\un\mu_\la$ distribution. The process $\un\ze(\cdot)$ initially 
has distribution $\un\mu_\vr$, except at sites $\{-n+1,\,\dots,\,0\}$ 
where the parameter $\vr$ has been replaced by $\la$.

Define a third initial configuration by 
\[
\xi_i(0)=\left\{
\ba
&\ze_i(0)&&\text{if }i\leq-n,\\
&\eta_i(0)&&\text{if }i>-n.
\ea
\right.
\]
We let all these processes evolve jointly in the basic 
coupling. The following majorizations are true initially 
and are preserved by the evolution:
\[
\un\eta(t)\leq\un\xi(t)\leq\un\ze(t)\qquad\text{and}\qquad \un h^\ze(t)\leq\un h^\xi(t)
\]
where the last inequality is for  column heights. 

Let us denote the net particle currents \eqref{eq:jvdef} by $J^{V,\eta}$ and 
$J^{V,\ze}$ in the respective processes $\un\eta(\cdot)$
and $\un\ze(\cdot)$. The first observation is that 
 $Q^{(-n)}$ gives one-sided control over the difference 
of  these currents. 
\begin{lm}\label{lm:geom}
There is a nonnegative  process $N(t)$ 
with constant  Ge\-o\-met\-ric($q/p$) time marginals 
$\Pv[N(t)=k]=(1-q/p)(q/p)^k$  and  
such that for any $V\in\Rb$ 
\[
Q^{(-n)}(t)\leq[Vt] \quad\text{implies}\quad
J^{V,\ze}(t)-J^{V,\eta}(t)\leq N(t).
\]
\end{lm}
\begin{proof}
Denote the positions of the  $\xi-\eta$ second class particles at time
$t$  by
\[ \dotsm <Y_k(t)< \dotsm <Y_{-2}(t)<Y_{-1}(0)<Y_0(t).\]
As argued before, we can arrange for these particles to keep
their labels, and then the ordering is preserved.  Initially
$Y_0(0)\leq -n=Q^{(-n)}(0)$. 

Let $m_Q(t)=\max\{k: Y_k(t)\leq Q^{(-n)}(t)\}$ 
be the label of the $\xi-\eta$ particle at or closest
to the left of $Q^{(-n)}(t)$.
  Initially $m_Q(0)=0$.
Once  $Q^{(-n)}(t)\in\{Y_k(t)\}$, this containment property will hold
forever because the basic coupling
will preserve  the ordering $\un\eta^+\leq\un\xi$ if this is ever
established. 

We claim that $m_Q$ remains zero while $Q^{(-n)}$ is disjoint
from the second class particles $\{Y_k\}$.   This follows
 from showing that 
 there is no jump  which swaps the ordering 
 $Y_0< Q^{(-n)}$. 
If $Y_0=i$ and $Q^{(-n)}=i+1$, then we have $0=\eta_i=\eta_i^+<\xi_i=1$ 
and $0=\eta_{i+1}=\xi_{i+1}<\eta_{i+1}^+=1$. In this situation
 column $i$ 
of $\xi$ can increase, or the same column of $\eta^+$ can decrease. Neither
 of 
these steps can swap the positions of $Y_0$ and $Q^{(-n)}$, 
instead they make $Y_0=Q^{(-n)}$.

Once $Q^{(-n)}$ is riding on the $\{Y_k\}$ particles so
that actually
\[
Q^{(-n)}(t)=Y_{m_Q(t)}(t),
\]
its label $m_Q$ evolves in the basic coupling as follows. 
\begin{itemize}
\item 
Suppose $Q^{(-n)}=Y_{k}=i$ and $Y_{k+1}=i+1$. A
 right Poisson arrow $(i\to i+1)$ appears at rate $p$
and  increases  $m_Q$
from $k$ to $k+1$. 
\item Suppose $Q^{(-n)}=Y_{k}=i$ and $Y_{k-1}=i-1$. 
A left Poisson arrow $(i\to i-1)$ appears at rate $q$
and decreases $m_Q$ 
from $k$ to $k-1$.  
\end{itemize}
When $\xi-\eta$ particle $Y_{m_Q(t)}$ itself jumps, $Q^{(-n)}$
jumps with it. 

To get bounds on $Q^{(-n)}(t)$ we introduce a suitable 
steady-state object.   As we did with the $\ze$-marks 
 \eqref{eq:zedef} for the upper bound, we introduce a further 
classification among  the $\xi-\eta$ second class particles
 so that exactly one of them has
priority over all the others.  This special particle is marked
by the label $m(t)$.   The Poisson arrows move $m(t)$ 
exactly the same way as $m_Q(t)$ on the labels $\{-\infty<k\leq 0\}$:
\begin{itemize}
\item If $Y_{m(t)}(t)=i$ and $Y_{m(t)+1}(t)=i+1$ and there is a right Poisson
arrow $(i\to i+1)$  then $m(t)$ increases by $1$.
 This happens  at rate $p$. 
\item If $Y_{m(t)-1}(t)=i-1$ and $Y_{m(t)}(t)=i$ and there is a 
left Poisson arrow $(i\to i-1)$ then $m(t)$ decreases by $1$.
 This happens  at rate $q$. 
\end{itemize}
From these rates we see that $m(t)$ behaves like a birth and 
death chain on $\Zb_-$ whenever adjacency of $\xi-\eta$ second class
particles permits $m(t)$ to jump.   
Without obstruction this birth and death
chain would have reversible measure $\pi(k)=(1-q/p)(q/p)^{\abs{k}}$ for 
$k\leq 0$.  
We give $m(0)$ initial distribution $\pi$.
Then the argument given for Proposition  \ref{pr:fks} implies that 
for each fixed time $t$ we have $\Pv[m(t)=k]=\pi(k)$. 

We have arranged  $m(0)\leq m_Q(0)$  initially. 
Since $Q^{(-n)}$ cannot hop over $Y_0$ without joining it,
the  identical responses of $m$ and $m_Q$ 
 to the Poisson arrows implies that 
$m(t)\leq m_Q(t)$ for all time.   

To connect with currents, 
note that  $-m_Q(t)$ equals  the number of $\xi-\eta$ 
second class particles strictly to the  
right of $Q^{(-n)}(t)$ at time $t$.  Since these particles started off
 in $(-\infty,0]$ and they account for all the 
discrepancies between the processes $\un\xi$ and $\un\eta$,  
this number equals 
 the height difference $h^\xi_{Q^{(-n)}(t)}(t)-
h^\eta_{Q^{(-n)}(t)}(t)$. Thus 
 under $\{Q^{(-n)}(t)\leq[Vt]\}$ we have
\[
\ba
-m(t){\geq}-m_Q(t)&=h^\xi_{Q^{(-n)}(t)}(t)-h^\eta_{Q^{(-n)}(t)}(t)\\
&=h^\xi_{[Vt]}(t)-h^\eta_{[Vt]}(t)+\sum_{i=Q^{(-n)}(t)+1}^{[Vt]}(\xi_i(t)-\eta_i(t))\\
&\geq h^\xi_{[Vt]}(t)-h^\eta_{[Vt]}(t)\\
&\geq h^\ze_{[Vt]}(t)-h^\eta_{[Vt]}(t)=J^{V,\ze}(t)-J^{V,\eta}(t).
\ea
\]
To obtain the statement of the lemma take $N(t)=-m(t)$.
\end{proof}

Let $\un\om(\cdot)$ be a stationary ASEP 
started from the  
$\un\mu_\vr$ Bernoulli distribution. The next lemma gives a way to
compare the distributions of  $\un\ze$ and 
$\un\om$. 
\begin{lm}\label{lm:rn}
Denote by $\Pv^\om$ and $\Pv^\ze$ the probability of events that depend only on 
the respective processes $\un\om(\cdot)$ and $\un\ze(\cdot)$. Then
\[
\Pv^\ze(\cdot)\leq\Pv^\om(\cdot)^\frac12\cdot\exp\Bigl[\frac{n(\vr-\la)^2}
{2\vr(1-\vr)}\Bigr].
\]
\end{lm}
\begin{proof}
Let 
\[
Z=\sum_{i=-n+1}^0\ze_i(0). 
\]
 $Z$ has a Binomial($n, \la$) distribution 
 $\nu^\la(z)=\binom{n}{z}\la^z(1-\la)^{n-z}$
for $z=0,\dots, n$.  We use the Cauchy-Schwarz inequality below to perform 
a change of measure on this binomial distribution. 
The binomial mass functions 
in the second line are easily added up.
\[
\ba
\Pv^\ze(\cdot)&=\sum_{z=0}^n\Pv^\ze(\,\cdot\,|\,Z=z)\,
[\nu^\vr(z)]^\frac12\cdot\frac{\nu^\la(z)}{[\nu^\vr(z)]^\frac12}\\
&\leq\Bigl[\sum_{z=0}^n[\Pv^\ze(\,\cdot\,|\,Z=z)]^2\,
\nu^\vr(z)\Bigr]^\frac12\cdot\Bigl[\sum_{z=0}^n\frac{[\nu^\la(z)]^2}
{\nu^\vr(z)}\Bigr]^\frac12\\
&=\Bigl[\sum_{z=0}^n[\Pv^\ze(\,\cdot\,|\,Z=z)]^2\,\nu^\vr(z)\Bigr]^\frac12
\cdot\Bigl[1+\frac{(\vr-\la)^2}{\vr(1-\vr)}\Bigr]^\frac{n}{2}\\
&\leq\Bigl[\sum_{z=0}^n\Pv^\ze(\,\cdot\,|\,Z=z)\,\nu^\vr(z)\Bigr]^\frac12
\cdot\exp\Bigl[\frac{n(\vr-\la)^2}{2\vr(1-\vr)}\Bigr].
\ea
\]
All that is left is to recognize that $\Pv^\ze(\,\cdot\,|\,Z=z)$ is 
the probability of a process $\un\ze(\cdot)$ whose initial 
distribution 
is Bernoulli$(\vr)$ outside $\{-n+1\dots0\}$, 
with $z$ particles distributed 
in that interval with each configuration equally likely. Summing 
these with the Binomial$(n,\,\vr)$ coefficients $\nu^\vr(z)$ 
gives the Bernoulli $\un\mu_\vr$ initial distribution of the
 process $\un\om(\cdot)$.
\end{proof}


\hop
{\bf Proof of the lower bound}

\hip
With these preparations we are ready to prove the lower bound. 
As for the upper bound, $Q_a(t)$ is the position 
of a second class antiparticle 
started from the origin on a $\un\mu_\vr$-equilibrium process $\un\om(\cdot)$ 
 initially perturbed by setting $\om_0(0)=1$. 
Our target quantity is abbreviated 
as before by $\Psi(t)=\Ev(|Q_a(t)-[V^\vr t]|)$. We start with bounding the 
probability of the complement of the event in Lemma \ref{lm:geom}.
As before $u$ is an arbitrary but fixed positive integer
 and $n=[V^\la t]-[V^\vr t]+u$.   
\begin{lm}
\[
\Pv\{Q^{(-n)}(t)>[V^\vr t]\}\leq\frac{\Psi(t)}{u}+\frac{4t(p-q)(\vr-\la)+2}{u}
+\frac{C}{u(\vr-\la)}.
\]
\label{lm:LBlm3} \end{lm}
\begin{proof}
For this proof, $Q(t)$ will refer to a second class particle started from 
the origin on a process $\un{\wih\eta}(\cdot)$ in $\un\mu_\la$ distribution, 
except for $\wih\eta_0(0)=0$. 
Translation invariance implies
\begin{multline*}
\Pv\{Q^{(-n)}(t)>[V^\vr t]\}=\Pv\{Q(t)-[V^\la t]>u\}\leq
\frac{\Ev(|Q(t)-[V^\la t]|)}{u}\\
\leq\frac{\Ev(|Q(t)-Q_a(t)|)}{u}+\frac{\Ev(|Q_a(t)-[V^\vr t]|)}{u}+\frac{[V^\la t]-[V^\vr t]}{u}
\end{multline*}
where we introduced  the couplings of Section \ref{sc:ub}. 
As in Lemma \ref{lm:vardif}, 
the first term is bounded from above by $2t(p-q)(\vr-\la)/u+C/[u(\vr-\la)]$. 
The second term is $\Psi(t)/u$, and the third term is bounded from above by 
$2t(p-q)(\vr-\la)/u+2/u$.
\end{proof}
\begin{lm}
For any $0<K<(p-q)t(\vr-\la)^2$,
\[
\ba
\Pv\{Q^{(-n)}(t)\leq[V^\vr t]\}&\leq\frac{\vr^{1/2}(1-\vr)^{1/2}\Psi(t)^{1/2}}
{(p-q)t(\vr-\la)^2-K}\cdot\exp\Bigl[\frac{n(\vr-\la)^2}{2\vr(1-\vr)}\Bigr]\\
&\quad+\frac{4\la \Psi(t)}{K^2}+
\frac{8t(p-q)\la(\vr-\la)}{K^2}+\frac{C}{K^2(\vr-\la)}
+\frac{C}{K-8}.
\ea
\]
\label{lm:LBlm4} \end{lm}
\begin{proof}
Lemma \ref{lm:geom} leads to
\begin{align}
&\Pv\{Q^{(-n)}(t)\leq[V^\vr t]\}\leq\Pv\{J^{V^\vr,\ze}(t)-J^{V^\vr,\eta}(t)
\leq N(t)\}\notag\\
&\qquad
\leq\Pv\{J^{V^\vr,\ze}(t)\leq K+(p-q)t(2\vr\la-\la^2)-1\}\label{eq:torn}\\
&\qquad\qquad
+\Pv\Bigl\{J^{V^\vr,\eta}(t)-(p-q)t(2\vr\la-\la^2)>\frac{K}{2}+3\Bigr\}
\label{eq:tovarj}\\
&\qquad\qquad
+\Pv\Bigl\{N(t)>\frac{K}{2}-4\Bigr\}.\label{eq:togeo}
\end{align}

We apply Lemma \ref{lm:rn} to line \eqref{eq:torn} to bound it by the 
equilibrium $\un\mu_\vr$-pro\-ba\-bi\-li\-ty: 
\[
\ba
&\quad\ \bigl[\Pv^\vr\{J^\vr(t)\leq K+(p-q)t(2\vr\la-\la^2)-1\}
\bigr]^\frac12
\cdot\exp\Bigl[\frac{n(\vr-\la)^2}{2\vr(1-\vr)}\Bigr]\\
&\leq\bigl[\Pv^\vr\{\wt J^\vr(t)\leq K-(p-q)t(\vr-\la)^2\}\bigr]^\frac12
\cdot\exp\Bigl[\frac{n(\vr-\la)^2}{2\vr(1-\vr)}\Bigr]\\
&\leq\frac{[\Vv(J^\vr(t))]^{1/2}}{(p-q)t(\vr-\la)^2-K}
\cdot\exp\Bigl[\frac{n(\vr-\la)^2}{2\vr(1-\vr)}\Bigr]\\
&=\frac{\vr^{1/2}(1-\vr)^{1/2}\Psi(t)^{1/2}}{(p-q)t(\vr-\la)^2-K}
\cdot\exp\Bigl[\frac{n(\vr-\la)^2}{2\vr(1-\vr)}\Bigr].
\ea
\]
The term $-1$ was subsumed
in errors caused by integer parts when we centered $J^\vr(t)$.

A simple coupling consideration shows that $\Ev(J^{V^\vr,\eta}(t))$ 
differs by at most one from the same expectation taken under
an unperturbed $\un\mu_\la$ initial condition. Thus taking integer parts 
again into account, line \eqref{eq:tovarj} is bounded by
\[
\Pv\Bigl\{\wt J^{V^\vr,\eta}(t)>\frac{K}{2}\Bigr\}\leq
\frac{4\Vv(J^{V^\vr,\eta})}{K^2}.
\]
 Lemma \ref{lm:vardif} can be applied to  bound this variance
even though the  second class particle now starts
 at $-n$ rather than at the origin.   The only change needed in the 
proof of Lemma \ref{lm:vardif} is in the calculation \eqref{eq:temp12}
where one must condition on $\eta_{-n}(0)$ instead of on $\eta_{0}(0)$. 
Hence we can continue from above to bound 
 line \eqref{eq:tovarj} with 
\[
\frac{4\la \Psi(t)}{K^2}+\frac{8t(p-q)\la(\vr-\la)}{K^2}+\frac{C}{K^2(\vr-\la)}.
\]

Lastly, the geometric probability \eqref{eq:togeo} is bounded by 
$C/(K-8)$ by Chebyshev's inequality.
\end{proof}

Now the last step of   the 
 lower bound of Theorem \ref{tm:main}.  By H\"older's
or Jensen's inequality it suffices to prove the case $\mom=1$,
in other words that
\[
\liminf_{t\to\infty} t^{-2/3}\Psi(t)>0.
\]
In the two last lemmas take 
\[
u=\lc ht^{2/3}\rc,\quad\vr-\la=bt^{-1/3}, \quad \text{ and }\quad K=bt^{1/3},
\]
where $h$ and  $b$ are large, in particular $b$ large enough 
to have  $b<(p-q)b^2$ so that $K$ satisfies Lemma \ref{lm:LBlm4}. 
Then
\[n=[V^\la t]-[V^\vr t]+u\leq[2(p-q)b+h]t^{2/3}+3
\le Ct^{2/3}\]
for large enough $t$. 
We can simplify the outcomes of Lemma \ref{lm:LBlm3} and
Lemma \ref{lm:LBlm4}  to the inequalities 
\be
\Pv\{Q^{(-n)}(t)>[V^\vr t]\}\leq C\frac{\Psi(t)}{t^{2/3}}
+\frac{4b}{h} +\frac{C}{t^{1/3}}
\label{eq:LBaux3}\ee
and
\be\begin{split} 
\Pv\{Q^{(-n)}(t)\leq[V^\vr t]\}&\leq
C\biggl(\frac{\Psi(t)}{t^{2/3}}\biggr)^{1/2}
+C\frac{\Psi(t)}{t^{2/3}}+
\frac{8}{b}+\frac{C}{t^{1/3}}.
\end{split}
\label{eq:LBaux4}\ee
The new constant $C$ depends on $b$ and $h$. 

The lower bound now follows because the left-hand sides of 
\eqref{eq:LBaux3}--\eqref{eq:LBaux4} add up to one for each 
fixed $t$, while we can fix   
 $b$ large enough and then $h$ large enough so 
that $4b/h + 8/b < 1$.  Then $t^{-2/3}\Psi(t)$ must
have a positive lower bound for all large enough $t$.
This completes the proof of
Theorem \ref{tm:main}. 

\smallskip

{\sl Remark.} The key properties, namely $f''\ne 0$
for the flux and asymmetry of the jump kernel, 
 were used somewhat surreptitiously in the above proof.
As in   calculation \eqref{eq:H''point}, 
 the nonvanishing of $f''$  renders 
$\Ev[J^\vr(t)-J^{V^\vr,\eta}(t)]$ of order $t^{1/3}$.  
This  allowed us to  take $K$ of order 
$t^{1/3}$,  then $\vr-\la$ of order $t^{-1/3}$, and 
 $u$ of order $t^{2/3}$.  This way  the final bounds 
 \eqref{eq:LBaux3}--\eqref{eq:LBaux4}
have  $\Psi(t)$ divided by the correct order $t^{2/3}$.

   Asymmetry also came in: if $p=q=1/2$ then 
the means of the equilibrium currents are zero  and there
would be no deviations to take advantage of
in the proof of  Lemma \eqref{lm:LBlm4}. 

\hop
\begin{center}
\footnotesize{
\begin{tabular}{l}
{\sc M\'arton Bal\'azs}\\
{\sc MTA-BME Stochastics Research Group, Institute of Mathematics}\\
{\sc Budapest University of Technology and Economics}\\
{\sc 1 Egry J\'ozsef u., $5^{\text{th}}$ floor 7, Bld.\ H}\\
{\sc 1111 Budapest, Hungary}\\
{\tt balazs@math.bme.hu}\\
\\
{\sc Timo Sepp\"al\"ainen}\\
{\sc Mathematics Department, University of Wisconsin-Madison}\\
{\sc Van Vleck Hall, 480 Lincoln Dr, Madison WI 53706-1388, USA.}\\
{\tt seppalai@math.wisc.edu}
\end{tabular}
}
\end{center}

\bibliography{refsmarton}

\begin{thebibliography}{10}

\bibitem{aldo-diac}
D.~Aldous and P.~Diaconis.
\newblock Hammersley's interacting particle process and longest increasing
  subsequences.
\newblock {\em Probab. Theory Related Fields}, 103(2):199--213, 1995.

\bibitem{bdj}
J.~Baik, P.~Deift, and K.~Johansson.
\newblock On the distribution of the length of the longest increasing
  subsequence of random permutations.
\newblock {\em J. Amer. Math. Soc.}, 12(4):1119--1178, 1999.

\bibitem{fluct}
M.~Bal{\'a}zs.
\newblock Growth fluctuations in a class of deposition models.
\newblock {\em Ann. Inst. H. Poincar{\'e} Probab. Statist.}, 39:639--685, 2003.

\bibitem{third}
M.~Bal{\'a}zs, E.\ Cator, and T.\ Sepp{\"a}l{\"a}inen.
\newblock Cube root fluctuations for the corner growth model associated to the
  exclusion process.
\newblock {\em Electronic Journal of Probability}, 11:1094--1132, 2006.

\bibitem{raprwre}
M.~Bal{\'a}zs, F.~Rassoul-Agha, and T.~Sepp{\"a}l{\"a}inen.
\newblock The random average process and random walk in a space-time random
  environment in one dimension.
\newblock {\em Comm. Math. Phys.}, 266(2):499--545, 2006.

\bibitem{varj2nd}
M.~Bal{\'a}zs and T.~Sepp{\"a}l{\"a}inen.
\newblock Exact connections between current fluctuations and the second class
  particle in a class of deposition models.
\newblock {\em Journal of Stat. Phys.}, 127(2):431--455, 2007.

\bibitem{cuberoot}
E.~Cator and P.~Groeneboom.
\newblock Second class particles and cube root asymptotics for {H}ammersley's
  process.
\newblock {\em Ann. Probab.}, 34(4), 2006.

\bibitem{dgl}
D.~D{\"u}rr, S.~Goldstein, and J.~Lebowitz.
\newblock Asymptotics of particle trajectories in infinite one-dimensional
  systems with collisions.
\newblock {\em Comm. Pure Appl. Math.}, 38(5):573--597, 1985.

\bibitem{evan}
L.~C. Evans.
\newblock {\em Partial differential equations}, volume~19 of {\em Graduate
  Studies in Mathematics}.
\newblock American Mathematical Society, Providence, RI, 1998.

\bibitem{se}
P.~A. Ferrari and L.~R.~G. Fontes.
\newblock Current fluctuations for the asymmetric simple exclusion process.
\newblock {\em Ann. Probab.}, 22:820--832, 1994.

\bibitem{fks}
P.~A. Ferrari, C.~Kipnis, and E.~Saada.
\newblock Microscopic structure of travelling waves in the asymmetric simple
  exclusion process.
\newblock {\em Ann. Probab.}, 19(1):226--244, 1991.

\bibitem{ferspohn}
P.~L. Ferrari and H.~Spohn.
\newblock Scaling limit for the space-time covariance of the stationary totally
  asymmetric simple exclusion process.
\newblock {\em Comm. Math. Phys.}, 265(1):1--44, 2006.

\bibitem{1/3}
K.~Johansson.
\newblock Shape fluctuations and random matrices.
\newblock {\em Comm. Math. Phys.}, 209:437--476, 2000.

\bibitem{cl}
C.~Kipnis and C.~Landim.
\newblock {\em Scaling limits of interacting particle systems}.
\newblock Springer-Verlag, Berlin, 1999.

\bibitem{lqsy}
C.~Landim, J.~Quastel, M.~Salmhofer, and H.-T. Yau.
\newblock Superdiffusivity of asymmetric exclusion process in dimensions one
  and two.
\newblock {\em Comm. Math. Phys.}, 244(3):455--481, 2004.

\bibitem{ips}
T.~M. Liggett.
\newblock {\em Interacting particle systems}.
\newblock Springer-Verlag, 1985.

\bibitem{stochi}
T.~M. Liggett.
\newblock {\em Stochastic interacting systems: contact, voter and exclusion
  processes}.
\newblock Springer-Verlag, 1999.

\bibitem{spohn}
M.~Pr{\"a}hofer and H.~Spohn.
\newblock Current fluctuations for the totally asymmetric simple exclusion
  process.
\newblock In {\em In and out of equilibrium (Mambucaba, 2000)}, volume~51 of
  {\em Progr. Probab.}, pages 185--204. Birkh{\"a}user Boston, Boston, MA,
  2002.

\bibitem{quava2}
J.~Quastel and B.~Valk{\'o}.
\newblock A note on the diffusivity of finite-range asymmetric exclusion
  processes on {$\mathbb Z$}.
\newblock {\em {\tt http://arxiv.org/abs/0705.2416}}, 2007.

\bibitem{quava}
Jeremy Quastel and Benedek Valko.
\newblock {$t\sp {1/3}$} {S}uperdiffusivity of finite-range asymmetric
  exclusion processes on {$\mathbb Z$}.
\newblock {\em Comm. Math. Phys.}, 273(2):379--394, 2007.

\bibitem{hl}
F.~Rezakhanlou.
\newblock Hydrodynamic limit for attractive particle systems on {$\Zb^d$}.
\newblock {\em Comm. Math. Phys.}, 140(3):417--448, 1991.

\bibitem{sepp-ejp}
T.~Sepp{\"a}l{\"a}inen.
\newblock A microscopic model for the {B}urgers equation and longest increasing
  subsequences.
\newblock {\em Electron. J. Probab.}, 1:no.\ 5, approx.\ 51 pp.\ (electronic),
  1996.

\bibitem{hkl}
T.~Sepp{\"a}l{\"a}inen.
\newblock Existence of hydrodynamics for the totally asymmetric simple
  {$K$}-exclusion process.
\newblock {\em Ann. Probab.}, 27(1):361--415, 1999.

\bibitem{flucha}
T.~Sepp{\"a}l{\"a}inen.
\newblock Second-order fluctuations and current across characteristic for a
  one-dimensional growth model of independent random walks.
\newblock {\em Ann. Probab.}, 33:759--797, 2005.

\bibitem{spi}
F.~Spitzer.
\newblock Interaction of {M}arkov processes.
\newblock {\em Adv. in Math.}, 5:246--290, 1970.

\bibitem{vBKS85}
H.~van Beijeren, R.~Kutner, and H.~Spohn.
\newblock Excess noise for driven diffusive systems.
\newblock {\em Phys. Rev. Lett.}, 54(18):2026--2029, May 1985.

\bibitem{yau}
H.-T. Yau.
\newblock {$(\log t)\sp {2/3}$} law of the two dimensional asymmetric simple
  exclusion process.
\newblock {\em Ann. of Math. (2)}, 159(1):377--405, 2004.

\end{thebibliography}
\bibliographystyle{plain}

\end{document}